\documentclass[preprint,12pt]{elsarticle}

 \usepackage[T1]{fontenc}
\usepackage{amsmath}
\usepackage{graphicx}
\usepackage{mathtools}   
\usepackage{amssymb}

\usepackage{stmaryrd}
\usepackage{bm}
\usepackage{amsfonts} 
\usepackage{float}
\usepackage{xcolor}
\usepackage{amsmath,bm}
\usepackage{mathabx}
\newcommand*\xbar[1]{%
  \hbox{%
    \vbox{%
      \hrule height 0.5pt 
      \kern0.5ex
      \hbox{%
        \kern-0.1em
        \ensuremath{#1}%
        \kern-0.1em%
      }%
    }%
  }%
} 
\newtheorem{theorem}{Theorem}[section]
\newtheorem{lemma}[theorem]{Lemma}
\newtheorem{proposition}[theorem]{Proposition}
\newtheorem{corollary}[theorem]{Corollary}

\newtheorem{definition}[theorem]{Definition}

\newenvironment{proof}[1][Proof]{\begin{trivlist}
\item[\hskip \labelsep {\bfseries #1}]}{\end{trivlist}}

\newenvironment{example}[1][Example]{\begin{trivlist}
\item[\hskip \labelsep {\bfseries #1}]}{\end{trivlist}}
\newenvironment{remark}[1][Remark]{\begin{trivlist}
\item[\hskip \labelsep {\bfseries #1}]}{\end{trivlist}}

\usepackage{amssymb}
\usepackage[nodots]{numcompress}
\makeatletter
\newcommand{\ostar}{\mathbin{\mathpalette\make@circled\star}}
\newcommand{\make@circled}[2]{%
  \ooalign{$\m@th#1\smallbigcirc{#1}$\cr\hidewidth$\m@th#1#2$\hidewidth\cr}%
}
\newcommand{\smallbigcirc}[1]{%
  \vcenter{\hbox{\scalebox{0.77778}{$\m@th#1\bigcirc$}}}%
}
\makeatother

\begin{document}

\begin{frontmatter}

\tnotetext[label1]{Corresponding author: Terence R. Smith, Department of Computer Science, University of California at Santa Barbara, Ca 93106}
\author{Terence R. Smith { }\fnref{label2}}
 \ead{smithtr@cs.ucsb.edu}
\title{Rational Wave Numbers and the Algebraic Structure of the Cyclic Groups of the Roots of Unity}
\begin{abstract}
Rational wave numbers are periodic sequences 
${\bm \omega}={\bf A}{\bf w}(f,g)$ in which
amplitude ${\bf A}$ a product of
powers of trigonometric sequences
and ${\bf w}(f,g)=\exp({\bf {i2}\pi ( f {\bm \xi} \oplus g{\bf 1})})$ 
is a sequence with
$  \xi \epsilon \mathbb{ Z}$ and
$f,g$ rational.
They generalize the cyclic groups of the $n$th roots of unity
and are generated from two unitary sequences.
The multiplicative group ${\bf W}_M$ 
with ${\bf A}={\bf 1}$
is their closure 
wrt
product, root, and reflection operators.
The commutative ring ${\bf W}_A$  
has additional closure 
wrt summation. The field
${\bf W}_I$ of invertible wave numbers has
further closure wrt an inverse.
Sums and differences of its elements are
$
{\bm \omega}_1 \oplus {\bm \omega}_2= {\bf 2}  \cos\big( {{\bf i}}
 ln\big(\frac{{\bm \omega}_2}{{\bm \omega}_1}\big)^{1/2}\big)
\big({\bm \omega}_1  {\bm \omega}_2 \big)^{1/2}$
and 
$
{\bm \omega}_1 \ominus {\bm \omega}_2= {\bf {2 i}}  \sin\big( {{\bf i}}
 ln\big(\frac{{\bm \omega}_2}{{\bm \omega}_1}\big)^{1/2}\big)
\big({\bm \omega}_1 {\bm \omega}_2 \big)^{1/2}$.
Its amplitudes 
form a multiplicative subgroup 
over which ${\bf W}_M \cup \big\{{\bf 0}\big\}$ is a vector space. 
Wave numbers of period $n\ \epsilon\ \mathbb{N}$ possess $n$ phases, multiplicative norms, 
representations wrt orthonormal bases, and prime representations.
Rational wave numbers may be completed
with respect to Cauchy sequences of parameters $({{f,g}})$.
Equations 
in invertible wave numbers
have solutions corresponding to zeros of their trigonometric factors.
Orthonormal bases are
employed in constructing 
the integral wave numbers and
allow definitions of the particulate and prime wave numbers
in terms of permissible phase-values.
\end{abstract}

 
\end{frontmatter}

\newpage

\section{Introduction}

The $n$th roots of unity arise in numerous branches
of mathematics 
(see, for example, Lang (1992\nocite{L92})).
It is well-known that each of the sets of the $n$th roots
of unity has the structure of a cyclic group, 
but less well-known that 
they may be generated from a set of two unitary elements.

It appears to be even less well-known that the unitary elements also
generate a nested sequence comprising a
multiplicative group, an additive group, and a field, all of whose elements 
have the polar form
${\bm \omega}={\bf A}{\bf w}(f,g)$ 
in which ${\bf w}(f,g) =e^{2\pi ( f {\bm \xi} \oplus g{\bf 1})}$ for
$ \xi \ \epsilon\ \mathbb Z; f,g \ \epsilon\ \mathbb Q$.
The amplitudes ${\bf A}$ of these elements are products of
functions of complex-valued trigonometric sequences
and are generated by the addition
of multiplicative wave numbers ${\bf w}(f,g)$.

Proofs of these assertions are based on representations
of wave numbers
as periodic sequences of complex numbers.
Hence the next section of the paper provides a brief summary 
of periodic sequences of complex numbers and their operators.

The third section characterizes the 
Abelian group of multiplicative wave numbers 
as the closure of the set of two unitary wave numbers under a set of 
four operators that includes the sequence product operator $\otimes$. 
The elements ${\bf w}(f,g)$ 
of this group
may be viewed in terms of the cyclic groups of the $n$th roots of unity
under permissible rotations and phase shifts
for $n\ \epsilon \ \mathbb{N}$.

The fourth section describes the commutative ring of additive wave
numbers that is formed by adding the sequence summation operator
$\oplus$ to the set of operators that generate the multiplicative wave
numbers.
A key result of this section is that an arbitrary sum
of additive wave numbers has a unique representation 
in the polar form ${\bf A}{\bf w}$ in which ${\bf A}$ is a product 
of functions of trigonometric sequences
and ${\bf w}$ is a multiplicative wave number.

The fifth section extends the set of additive wave numbers to the set
of invertible wave numbers with the addition of an inverse operator.
It is shown that the set of invertible wave numbers forms a field 
while the set
of their amplitudes forms a multiplicative subgroup of the field.
The field may be represented as a vector space formed by
the cyclic groups of $n$th roots of unity 
over the elements of the group of amplitudes.

The product form of the amplidudes of this field
leads naturally to theorems concerning 
the number of solutions to equations in invertible wave numbers.
It is further shown how one may construct
an orthonormal basis for the wave numbers
and illustrates its application in the construction of classes of wave numbers
that include the integral, particulate,
and prime wave numbers.
Invertible wave numbers are characterized 
by a period $n\ \epsilon \ \mathbb{N}$
and $n$ phases unless constrained by co-numbers 
to have the form of the aperiodic particulate wave numbers. 

Interpreting wave numbers and their arithmetic
as an elementary theory of waves suggests
basing a theory of numbers
on the concept of a wave rather than the concept of a set
(as, for example, in Goldrei (1996)), as well as suggesting
applications relating to the wave equations of physics.
\section{Periodic Sequences of Complex Numbers
and their Operators
}\label{sec1}
\label{sec:seq}
\begin{definition}
\label{def:linphasenum}
An exponential function of a phase index
${\bm \rho}(f,g)=f {\bm \xi} \oplus g {\bm 1}$,
in which ${\bm \xi}$ is the sequence of integers $\mathbb{Z}$,
${\bf 1}$ is a
countable sequence of $1$'s, and $f,g$ are numbers,
is the sequence of exponential terms
\begin{eqnarray}
e^{2 \pi i( f {\bm \xi}+ g{\bf 1} )}=
\big\{ e^{    2 \pi i ( f {\xi}+ g)}
\mid \xi\ \epsilon\ \mathbb{Z} \big\},
\end{eqnarray}
with the function said to be rational if $f,g\ \epsilon \ \mathbb{Q}$.

\end{definition}
The equality of two exponential sequences
is defined in an element-wise manner and
operators and functions are applied
in an element-wise manner.
Representing a general sequence as 
$
{{\bf C}}= \Big\{{c}(\xi) \Big\}
$
the admissible operators are
\begin{enumerate}
\item
conjugate sequence operator:  
$\bar{{\bf C}}=
{\it M_R}({\bf C}) = \Big\{\bar{c}(\xi) \Big\}
$;
\item
orthogonal conjugate sequence operator: 
$
\tilde{{\bf C}}={\it M_I}({\bf C}) =\Big\{-\bar{c}(\xi) \Big\}
$;
\item
$n$th root operator: 
$
{\it R_{n}}({\bf C})
=\Big\{ c(\xi)^{1/n} \Big\}, \hspace{.2in} n\ \epsilon \ \mathbb{N};
$
\item
product operator $\otimes$
applied to pairs of elements
with identical  $\xi$-values:
$
{\bf C}_1 \otimes{\bf C}_2= \Big\{ c_1(\xi) \times c_2(\xi)\Big\}  
$;
\item
sum ${\bf C}_1 \oplus {\bf C}_2$
operator applied to pairs of elements
with identical $\xi$-values:
${\bf C}_1 \oplus {\bf C}_2 
=\Big\{ c_1(\xi) + c_2(\xi)\Big\}$;
\item
inverse operator:  
$
{\it I}\big({\bf C}\big) ={\it I}\Big(\Big\{ c(\xi) \Big\}\Big)
=\Big\{ \tfrac{1}{c(\xi)} \Big\},
$
which is not permitted on any sequence that contains
one or more zero elements.
\end{enumerate}
The product symbol $\otimes$
is omitted when there is no unambiguity.
One may define operators for
difference $({\bf C}_1\ominus {\bf C}_2) \coloneq {\bf C}_1 \oplus (-{\bf C}_2)$
and division
$({\bf C}_1\odiv {\bf C}_2)  \coloneq {\bf C}_1 \otimes {\it I}{(\bf C}_2)
$ which by the definition of the inverse operator disallows division by 
any sequence with one or more zero elements.
It is convenient to define the operator ${R}_{\infty}$ 
that computes
the root for every $n\ \epsilon \ \mathbb{N}$
\begin{eqnarray}
{R}_{\infty}({\bf C})=\big\{R_{n}({\bf C}) \mid 
n\ \epsilon \ \mathbb{N} \ \big\}
\end{eqnarray}
The commutative, associative, and distributive properties of the operators
follow immediately from the same properties of the complex numbers.
Operators may be composed and the compositions
used in defining additional operators, such as
the $n$th power of a sequence 
${\bf C}^n\coloneq   \overset {n} {\underset {k=1} \otimes } {\bf C} $
and the
conjugate product $\overline{\otimes}$ 
of two sequences $
{\bf C}_1\overline{\otimes}{\bf C}_2 \coloneq
{\bf C}_1 {\otimes}\overline{\bf C}_2$
whose definition implies
\begin{proposition}
\label{prop:inverse}
${\it I}({\bf C})=\overline{\bf C}\otimes{\it I}({\bf C}\overline{\otimes}{\bf C}) $
is the inverse  of a sequence ${\bf C}.$
\end{proposition}
A key property of the sequences discussed in this paper is their periodicity.
\begin{definition}
A sequence is periodic in ${\bm \xi}$ iff there is a finite integer $\xi_0$
s.t.
\begin{eqnarray}
C({\bm \rho}(f,g))=C({\bm \rho}(f,f\xi_0+g)),
\ \ \ \forall\ \xi\ \epsilon\ {\mathbb Z}.
\end{eqnarray}
\end{definition}
\begin{proposition}
\label{theor:period}
The product $\otimes$ and sum $\oplus$ of N sequences ${\bf C}_i, i=1,N$
with periodicities $d_1,...,d_N$ are periodic
with period
$p_N=\underset {j=1} {\overset{ N} {   \Pi} }d_j\Big/
\underset {j=1} {\overset{ N} {   \Pi} } g_j  $,
in which $g_1=1$ and 
$g_j =gcd\Big(\underset {l=1} {\overset{ j-1} {   \Pi} }d_l\Big/
\underset {l=1} {\overset{ j-1} {   \Pi} g_l}, d_j\Big)$
and $gcd(n,m)$ is the greatest common divisor of $n$ and $m$. 
\begin{proof}
Considering first the case of two periodic sequences ${\bf C}_1, {\bf C}_2$,
one notes that their
$\xi th$ elements satisfy the relations
$c_1(\xi)=c_1(\xi+d_1)$ and $c_2(\xi)=c_2(\xi+d_2)$
in which $d_1,d_2$ are their periods.
If $d_1$ and $d_2$ have a greatest common divisor $f=gcd(d_1,d_2)$,
then $d_1=f {d_1}^\prime$ and $d_2=f {d_2}^\prime$ 
for some ${d_1}^\prime < d_1$ and  ${d_2}^\prime < d_2$
and the sum of the $\xi$th elements of the two sequences
satisfies the relation
\begin{eqnarray}
c_1(\xi) +c_2(\xi) =c_1(\xi+n_1 f {d_1}^\prime) +c_2(\xi+n_2 f {d_2}^\prime)
\end{eqnarray}
for integers $n_1,n_2$ and it follows that the sum is periodic. 
Choosing $n_1=d_2/f$ and $n_2=d_1/f$, 
one obtains
\begin{eqnarray}
c_1(\xi) +c_2(\xi) 
 &=&c_1(\xi+{d_2}{d_1}^\prime) 
+c_2(\xi+{d_2}^\prime{d_1})
\\
 &=&c_1(\xi+f{d_2}^\prime{d_1}^\prime) 
+c_2(\xi+f{d_1}^\prime{d_2}^\prime)
\end{eqnarray}
and since ${d_1}^\prime$ and $ {d_2}^\prime$ are relatively prime,
$f{d_1}^\prime {d_2}^\prime=d_1d_2/f$ is the smallest period of the
sum ${\bf C}_1\oplus {\bf C}_2$. 
An analogous argument applies to the product
of the two sequences.
It follows that the value of the function
\begin{eqnarray}
\rho(p_i, p_j)= \frac{p_i p_j}{gcd (p_i,p_j)}
\end{eqnarray}
is the periodicity of the product or sum of the sequences
${\bf C}_i, {\bf C}_j$.\\
In generalizing this result to a set of $N$ sequences $\Big\{{\bf C}_1,...{\bf C}_N\Big\}$,
it is sufficient to apply $\rho(p_i,p_j)$ recursively to
the periodicities of the $N$ sequences
\begin{eqnarray}
\rho(d_1,d_2), \ \rho(\rho(d_1,d_2), d_3),\  \rho( \rho(\rho(d_1,d_2), d_3), d_4),...
\end{eqnarray}
and the results of the proposition follow immediately from a recursive
application of the function $\rho$ to the set
$\Big\{ {\bf C}_1, {\bf C}_2, ..., {\bf C}_N\Big\}$.
For the first three values of $N$, for example, one has
\begin{eqnarray}
p_1&=&n_1/g_1;\  g_1=1; \\
p_2&=&d_1d_2/g_1g_2; \ g_2=gcd(d_1/g_1, d_2);\nonumber\\
p_3&=&d_1d_2d_3/g_1g_2g_3;\  g_3=gcd(d_1d_2/g_1g_2,d_3)\nonumber 
\end{eqnarray}\qed 
\end{proof}
\end{proposition}
\begin{definition}
The principal subsequence of a periodic
sequence ${\bf C}$ 
is 
\begin{eqnarray}
\label{eq:def}
{_1} {\bf C}=  
\Big\{ c(\xi)\Big\}_{\xi=1}^{\xi=\lambda} 
\end{eqnarray}
in which $\lambda$ is the period of ${\bf C}$.
The $p$th principal subsequence 
${_p} {\bf C}$ is analogous.
\end{definition}
\begin{corollary}
\label{cor:sufficient}
Periodic sequences and 
the operations upon them
are finitely representable.
\begin{proof}  
A periodic sequence may be represented by its principle subsequence,
while the results of applying a unary operator
my be represented by the transformed principle sequence.
By Proposition \ref{theor:period}, the principle sequence 
of the sum and product 
of periodic sequences may be
represented by the sum and product of finite concatenations
of the principal sequences of the input sequences.
\end{proof}
\end{corollary}
The norm of a periodic sequence has
\begin{definition}
\label{def:norm}
The norm $\| {\bf C} \|$ of a sequence ${\bf C}$ is a mapping from the space of sequences
to the real numbers
\begin{eqnarray}	
\| {\bf C} \| 
=	\sqrt{     \frac{ \underset {\lambda}{\Sigma}		
               ( {\bf C}\overline{\bf C}) }	{\lambda}}
\end{eqnarray}	
in which $\lambda$ is the period of the sequence the sum is over
the phases of the principal sequence.
\end{definition} 
\begin{theorem} $  \| {\bf C}_1\otimes{\bf C}_2 \|  
=\|  {\bf C}_1\|\times \|{\bf C}_2\|$
\label{multiplicative}
%
\begin{proof}
\begin{eqnarray}
\label{eq:f122}
\lVert{\bf C}_1\rVert\times\lVert{\bf C}_2\rVert
=\sqrt{  \frac  { \underset {\lambda_{1}} {\Sigma}	c_1^2({ \xi})} {\lambda_1}  }
\times
\sqrt{  \frac  { \underset {\lambda_{2}} {\Sigma}	c_2^2({ \xi})} {\lambda_2}  }
=\sqrt{  \frac  { \underset {\lambda_{1}} {\Sigma}	(
c_1^2({ \xi}) \times  \underset {\lambda_{2}} {\Sigma}	c_2^2({ \xi}))} 
{\lambda_1  \lambda_2  } }.
\end{eqnarray} 
Noting that ${\bf C}_1$ and ${\bf C}_2$ are periodic sequences
with periods $\lambda_1$ and $\lambda_2$, 
the product wave number has period 
$\lambda_{12}=\lambda_{1}\lambda_{2}$, hence
\begin{eqnarray}
\label{eq:f222}
\lVert{\bf C}_1\otimes{\bf C}_2\rVert
=\sqrt{ \frac{ \underset {\lambda_{12}}{\Sigma}
( {\bf C}_1 {\bf C}_2\overline{ {\bf C}_1 {\bf C}_2}}{\lambda_{12}}}
=\sqrt{ \frac{ \underset {\lambda_{12}}{\Sigma}	(	{\bf C}_1\overline{\bf C}_1
\otimes{\bf C}_2\overline{\bf C}_2)}{\lambda_{12}}}
=\sqrt{ \frac{ \underset {\lambda_{12}}{\Sigma}	(	c_1^2({ \xi})	\times
	c_2^2({ \xi})	)}{\lambda_{12}}}	
\end{eqnarray}
in which the second equation follows from the first by commutativity.
 It follows from the definition
of the product $\otimes$ of sequences respectively of length 
$\lambda_{1}$ and $\lambda_{2}$ 
that the product of the sums in Equation \ref{eq:f122} has length $\lambda_{12}$
and 
is a factorization of the sum of the same length in Equation \ref{eq:f222}.
Hence the two expressions in Equations \ref{eq:f122} and \ref{eq:f222}
are equal.\qed
\end{proof}  
\end{theorem}
It is convenient to define unary and binary set operators
\begin{eqnarray}
\label{def:setop}
\widehat{\it O}_U({\bf S_1})&=&\big\{{\it O}_U(s)\ |\forall \ s\ \epsilon\ {\bf S}_1 \big\} \nonumber \\
\widehat{\it O}_B({\bf S}_1,{\bf S}_2)&=&\big\{{\it O}_B(s_1,s_2)\ |\forall \ s_1\ \epsilon\ {\bf S}_1, 
 \ s_2\ \epsilon\ {\bf S}_2 \big\}
\end{eqnarray}
that apply
to all members of a set of sequences $S$, as employed in
\begin{definition} $\ $ 
\label{def:multigen2}
The closure ${\bf S}$ of a finite set of sequences 
${\bf S}_{0}=\big\{{\bf C}_1,...,{\bf C}_m\big\}$ with respect
to a set operator 
$\widehat{O}$ is
\begin{eqnarray}
{\bf S} =	  {\underset {n \rightarrow \infty} {\lim}}  \  
 {\underset {j=0} {\overset {n} {\bigcup} }\widehat{O}^j({\bf S}_0)}
\end{eqnarray}
\end{definition}
One may obtain the following characterization of the generalization
of the closure with respect to several operators
\begin{proposition} $\ $ 
\label{def:multigen1}
The set of sequences 
$ {\bf S}=
{\bf S}_{0}\cup
 \Big(  \widehat{O}_1  ({\bf S}) \ \cup \ ...\cup
\widehat{ O}_n ({\bf S}) \Big) 
$
\ represents the closure of the set ${\bf S}_{0}$ with respect to the set operators
$\big\{\widehat{ O}_1,...,\widehat{ O}_n\big\}$.
\begin{proof}
Consider the terms $ {\bf S}=
{\bf S}_{0}\cup
 \widehat{O}_1  ({\bf S}) $ involving
the operator $ \widehat{O}_1 $. Since they imply that ${\bf S}_{0} \subset 
{\bf S}$ it follows that $\widehat{O}({\bf S}_{0}) \subset {\bf S}$ and, recursively
for all $n$, that  $\widehat{O}^n({\bf S}_{0}) \subset {\bf S}$.
Hence by Definition \ref{def:multigen2}
it follows that the closure of ${\bf S}_{0}$ with respect to $ \widehat{O}_1$ is contained in $ {\bf S}. $
Analogously, its closure with respect to the remaining operators
must also be contained in  $ {\bf S}. $
\qed
\end{proof}
\end{proposition}
This paper employs the set operators
$
\big\{   \widehat{\it M}_R,
{\widehat{\it M}}_I, 
{\widehat{\it R}}_N, 
{\widehat{\otimes}},
{\widehat{\oplus},
{\widehat{I}}}.
\big\}
$ corresponding to the sequence operators
$
\label{eq:ops}
\big\{ {\it M}_R,
{\it M}_I,
{\it R}_N,
\otimes, 
\oplus,
{\it I}.
\big\}
$
\section{The Abelian Group ${\bf W}_M(\otimes) $ of Multiplicative Wave Numbers}
In general, wave numbers are finitely-representable, 
periodic sequences of complex numbers
that may be generated from the
unitary wave numbers.
\begin{definition}
\label{def:multiplicative}
A wave number is unitary iff it is a member of the set 
\begin{eqnarray}
\label{eq:unitary}
{\bf W}_U
= \big\{{\bf w}(1,0)=
e^{2 \pi i {\bm \xi}},\ {\bf w}(0,1)=e^{2 \pi i {\bm 1}} \big\}.
\end{eqnarray}
in which 
${\bm \xi}$ is the sequence of all integers and
${\bf 1}$ is a countable sequence of $1$'s.
\end{definition}
\subsection{The Definition and Representation of Multiplicative Wave Numbers\\}
The multiplicative wave numbers ${\bf W}_M$
may be generated from
the set ${\bf W}_U$ of unitary wave numbers 
using the set of operators
${\bf S}_M= \big\{{M}_R, {M}_I, { R}_n,\otimes \big\}$.
\begin{definition} 
\label{def:multiplic}
A wave number ${\bm \omega}$
is multiplicative and ${\bm \omega}\ \epsilon\ {\bf W}_M $ iff
\begin{eqnarray}
{\bf W}_M ={\bf W}_U \cup
{\widehat{{M}}_R}({\bf W}_M)  \cup 
{\widehat{{M}}_I}({\bf W}_M)  \cup  {\widehat{{ R}}_N}({\bf W}_M)  \cup  
{\bf W}_M\widehat{\otimes}{\bf W}_M
\end{eqnarray}
in which 
${\widehat{{M}}_R},
{\widehat{{M}}_I}, {\widehat{{ R}}_N},
\widehat{\otimes}$ are set operators.
\end{definition}
It follows from Definition \ref{def:multiplic} that the
unitary wave numbers are multiplicative 
and it follows from Proposition \ref{def:multigen1}
that the multiplicative wave numbers ${\bf W}_M $ represent the closure
of ${\bf W}_U$
with respect to the operators
of ${\bf S}_M.$ 
Definition \ref{def:multiplic} implies that multiplicative wave numbers
may be represented as sequences of complex exponential terms: 
\begin{theorem}
\label{thm:ratform}
A wave number
is multiplicative if and only if it has a 
representation of the form
\begin{eqnarray}
\label{eq:ratform}
 {\bf w}(f,g)=  
e^{2 \pi i \big(f{\bm \xi}+ g{\bf 1}\big) },
\hspace{.2in} f,\ g\ \epsilon\ \mathbb{Q}.
\end{eqnarray}
\begin{proof} 
The closure of ${\bf W}_U $ with respect
to $\otimes$ implies that wave numbers of the form
$e^{2 \pi i {m}{\bm \xi} },
e^{2 \pi i p{\bf 1} },\ m,p\ \epsilon\ \mathbb{N}  $ are contained in ${\bf W}_M $;
closure with respect to $R_N$ implies that 
wave numbers of the form
$e^{2 \pi i \tfrac{m}{n}{\bm \xi} }$ and
$e^{2 \pi i \tfrac{p}{q}{\bf 1} }\ \tfrac{m}{n},\ \tfrac{p}{q}\ 
\epsilon\ \mathbb{Q}$ are contained in ${\bf W}_M $;
closure with respect to $M_I, M_R$ implies that 
wave numbers of the form
$\pm e^{\pm 2 \pi i \tfrac{m}{n}{\bm \xi} },
\pm e^{\pm 2 \pi i \tfrac{p}{q}{\bf 1} }  $ are contained in ${\bf W}_M $;
and closure with respect to $\otimes$ implies that 
wave numbers of the form
$\pm e^{2 \pi i \big(\pm \tfrac{m}{n}{\bm \xi} 
\pm \tfrac{p}{q}{\bf 1}\big) }  $ are contained in ${\bf W}_M $.
It follows from the properties of the exponential function 
and the definitions of the operators ${\bf S}_M$
that the form (\ref{eq:ratform}) is invariant under the operators since
\begin{eqnarray}
{M}_R( {\bf w}(f,g))&=&
{\bf w}(-f,-g) \\
{M}_I ( {\bf w}(f,g))&=&
{\bf w}(-f, \pi-g)   \nonumber \\
{ R}_N( {\bf w}(f,g))&=&
{\bf w}(\tfrac{f}{n},\tfrac{g}{n}) \nonumber  \\
 {\bf w}(f_1,g_1) \otimes{\bf w}(f_2,g_2)&=&
{\bf w}(f_1+f_2, g_1+g_2)  \nonumber
\end{eqnarray}
\end{proof}
\end{theorem}
\begin{corollary}
A multiplicative wave number 
${\bf w}(\frac{m}{n},\tfrac{p}{q})$
has a unique minimal representation. 
\begin{proof}
Applying the fundamental theorem of arithmetic 
to the multiplicative wave number 
${\bf w}(\frac{m}{n},\tfrac{p}{q})$
one obtains
the unique prime representations
$m={ \underset {i=1} {\overset{N_m} {   \Pi} }{p_i^{\alpha_i}}}$,
$n={ \underset {i=1} {\overset{N_n} {   \Pi} }{p_i^{\beta_i}}}$,
$p={ \underset {i=1} {\overset{N_p} {   \Pi} }{p_i^{\gamma_i}}}$,
$q={ \underset {i=1} {\overset{N_q} {   \Pi} }{p_i^{\delta_i}}}$,
in which $p_i$ is the ith prime number
and ${N_m},{N_n},{N_p},{N_q}$ represent the minimum prime number
indices required to represent
$m,n,p,q$. One has,
after cancellations of prime factors that leave only non-negative exponents,
\begin{eqnarray}
{\bf w}(\tfrac{m}{n},\tfrac{p}{q})=
{\bf w}
\Bigg(\frac{{ \underset {i=1} {\overset{N_m} {   \Pi} }{p_i^{\alpha_i}}}}
{ \underset {i=1} {\overset{N_n} {   \Pi} }{p_i^{\beta_i}}},
\frac{{ \underset {i=1} {\overset{N_p} {   \Pi} }{p_i^{\gamma_i}}}}
{ \underset {i=1} {\overset{N_q} {   \Pi} }{p_i^{\delta_i}}}\Bigg)
=
{\bf w}\Bigg(\frac{{ \underset {i=1} {\overset{N_m} {   \Pi} }{p_i^{\bar{\alpha}_i}}}}
{ \underset {i=1} {\overset{N_n} {   \Pi} }{p_i^{\bar{\beta}_i}}},\frac{{ \underset {i=1} {\overset{N_p} {   \Pi} }{p_i^{\bar{\gamma}_i}}}}
{ \underset {i=1} {\overset{N_q} {   \Pi} }{p_i^{\bar{\delta}_i}}}\Bigg)
\end{eqnarray} 
in which
$({\bar{\alpha}_i}\geq 0, {\bar{\beta}}_i=0)$ 
or $({\bar{\alpha}_i}= 0, {\bar{\beta}_i\geq0})$
and $({\bar{\alpha}_i}\geq 0, {\bar{\beta}_i=0})$
or   $({\bar{\alpha}_i}= 0, {\bar{\beta}_i\geq0})$
and 
$({\bar{\gamma}_i}\geq 0,         {\bar{\delta}_i=0})$   
or $({\bar{\gamma}_i}= 0,         {\bar{\delta}_i\geq0})$
and $({\bar{\gamma}_i}\geq 0, {\bar{\delta}_i=0})$ 
or $({\bar{\gamma}_i}= 0,         {\bar{\delta}_i\geq0})$. 
\end{proof} 
\end{corollary}
\begin{remark}
The numbers $(f,g)$, termed the parameters of a multiplicative wave number,
 are assumed to be always in minimal form and hence to have relatively prime 
numerators and denominators.
\end{remark}
\subsection{Elementary Properties of the Multiplicative Wave Numbers\\}
The elementary properties of multiplicative wave numbers
that are considered here
include their periodicity, primality, phase structure, group structure, 
rotational structure,
and multiplicative norm.
\begin{proposition}
A multiplicative wave number ${\bf w}(\tfrac{m}{n},\tfrac{p}{q})$ 
is a periodic function of period ${n}$
in the phase index ${\bm \xi}$.
\begin{proof}
\begin{eqnarray}
e^ {2 \pi i \big(\big(\tfrac{m}{n}\big)({\bm \xi}+n{\bf 1}) + \big(\tfrac{p}{q}\big){\bf 1}\big)}
= e^{2 \pi i\big(\tfrac{m}{n}{\bm \xi} + \tfrac{p}{q}{\bf 1}\big)} e^{2\pi i m{\bf 1}} 
= e^{2 \pi i \big(\tfrac{m}{n}{\bm \xi}+ \tfrac{p}{q}{\bf 1}\big)}.  
\end{eqnarray}
and $n$ is the minimum period since $\tfrac{m}{n}$ is in minimal form.
\end{proof} 
\end{proposition}
\begin{definition}
A multiplicative wave number 
${\bf w}(\frac{m}{n},\tfrac{p}{q})$ is period-prime iff its period $n$ is a prime number.
\end{definition}
\begin{theorem}
\label{equalityconst1}
A multiplicative wave number ${\bf w}(\frac{m}{n},\tfrac{p}{q}) $
may be represented as a product
of period-prime wave numbers
$
\Big(\underset {k=1} {\overset{n} {   \otimes} }
{\bf w}(\tfrac{m_k}{p_k},0)\Big)\Big(
\underset {k'=1} {\overset{n'} {   \otimes} }
{\bf w}(0,\tfrac{m^{'}_k}{p^{'}_k})\Big)$, 
unique up to order,
in which $n=\underset {k=1} {\overset{n} {   \Pi p_k} }$ and
$q=\underset {k=1} {\overset{n} {   \Pi q_k} }$ are prime representations.
\begin{proof}
Since the first and second arguments of ${\bf w}(\frac{m}{n},\tfrac{p}{q})$ 
may be converted independently into prime representations,
it suffices to show the conversion of the first argument $\tfrac{m}{n}$.
By the fundamental theorem of arithmetic,
let $n={ \underset {k=1} {\overset{N} {   \Pi} }{p_k}}$,
in which $p_k$ is the kth occurrence of a prime number and write
\begin{eqnarray}
\frac{m\ mod\ n}{n}=
\frac{m\ mod\ n}{{ \underset {k=1} {\overset{N} {   \Pi} }{p_k}}}
=\underset {k=1} {\overset{N} {   \Sigma} }\tfrac{m_k}{{p_k}}
=\tfrac{\underset {k=1} {\overset{N} {   \Sigma} }    
({m_k}     
{ \underset {j\ne k} {\overset{N} {   \Pi} }     {p_j )                    }          }}
{{ \underset {k=1} {\overset{N} {   \Pi} }         {p_k}}};
\end{eqnarray}
for some set of integers $m_k\ \epsilon \ \big\{1,...,p_k\big\}$
such that, given $m$, one has
\begin{eqnarray}
{m}
={\underset {k=1} {\overset{N} {   \Sigma} }    
({m_k}     
{ \underset {j\ne k} {\overset{N} {   \Pi} }     {p_j  )                    }          }}
{{      }}\ mod \ {{ \underset {k=1} {\overset{N} {   \Pi} }{p_k}}}
\end{eqnarray}
It follows that
\begin{eqnarray}
{\bf w}(\underset {k=1} {\overset{N} {   \Sigma} }\tfrac{m_k}{{p_k}},0)
= 
\underset {k=1} {\overset{N} {   \otimes} }{\bf w}(\tfrac{m_k}{{p_k}},0)
\end{eqnarray}
\qed
\end{proof}
\end{theorem}
\begin{definition}
The phases
of a multiplicative wave number ${\bf w}(\tfrac{m}{n}, \tfrac{p}{q}))$
are the $n$ equivalence classes represented by
$\Big\{e^{i2 \pi \big(\tfrac{( m { \xi}\ mod\  n) }{n}+ \tfrac{p}{q}\big)}\Big| \ \xi=1,n\Big\}$.
\end{definition}
\begin{proposition}
The phases of a multiplicative wave number and the elements of 
its principal sequence are isomorphic.
\begin{proof} 
\begin{eqnarray}
\Big\{e^{  i2 \pi \big(  \tfrac{( m { \xi}\ mod\  n) }{n}+ \tfrac{p}{q}\big)}
\Big| \ \  \xi=1,n \Big\}\ 
= \ 
\Big\{ e^  {  i 2 \pi \big( \tfrac{ m \xi}{n}+ \tfrac{p}{q} \big)   \big)}
\big|
 \ \  \xi=1,n  \Big\} 
\end{eqnarray}.
\end{proof} 
\end{proposition}
\begin{definition}
A wave number ${\bf w}(\tfrac{m}{n},g)$ is constant iff 
it has a single phase.
\end{definition}
\begin{theorem}
The multiplicative wave numbers 
have the structure of a multiplicative Abelian group ${\bf W}_M(\otimes)$
under
the product operator when they are represented in terms of their phases.
\begin{proof}
By Definition \ref{def:multiplic}
the set of multiplicative wave
numbers is closed under $\otimes$, which
is commutative and associative. 
The wave numbers ${\bf w} (m,n)
={\bf 1}
\ \forall\ m,n\ \epsilon\ {\mathbb Z}$
all represent the unique identity element 
${\bf 1}$ of ${\bf W}_M$ under $\otimes$
since
${\bf w} (m,n){\bf w}(f,g)= {\bf w}(m+f,n+g)
= {\bf w}(f,g)
$ and 
$e^{i2 \pi \big({( m {\bm \xi} }+ {n} {\bf 1} \big)}={\bf 1}$
is a unique phase.
The conjugate multiplicative wave number 
${\bf w}(-f,-g)$ 
is the unique multiplicative inverse of the wave number
${\bf w}(f,g)$, and
${\bf w}(f,g)\overline{\bf w}(f,g)
= {\bf 1}$, since 
$e^{i2 \pi \big(\tfrac{( m {\bm \xi}\ mod\  n) }{n}+ \tfrac{p}{q}{\bf 1}\big)}\otimes
e^{-i2 \pi \big(\tfrac{( m {\bm \xi}\ mod\  n) }{n}+ \tfrac{p}{q}{\bf 1}\big)}={\bf 1}$
\end{proof}
\end{theorem}
One may construct many subgroups of ${\bf W}_M(\otimes)$
by placing constraints on their generators or their generating operators, as
now exemplified in terms of 
the cyclic groups of the roots of unity, the simple multiplicative
wave numbers, and the constant multiplicative wave numbers.
\begin{proposition}
The cyclic group of the roots of unity for $n \ \epsilon\ \mathbb{Z}$
is the multiplicative wave number ${\bf w}(\tfrac{1}{n},0)$ generated from the phase
$e^{i2\pi \tfrac{{ 1}}{n}}$
by the set of operators $\big\{ M_R, \otimes \big\}$.
\end{proposition}
\begin{example}
\label{prop:i}
For $n=4$, one has the principal sequence
$
_{1}{  \bf w}(\tfrac{1}{4},0) 
=\big\{{i},\  -{1},\ 
 -{ i},\  
{ 1}\big\}.
$
\end{example}
\begin{proposition}
\label{lem:const}
The set ${\bf W}_{MC}$ of constant wave numbers forms a subgroup
of the multiplicative wave numbers
\begin{proof}
The set of constant numbers is closed under $\otimes$, has a unique
inverse element, and the  identity wave number ${\bf 1}$ 
has a single phase and is hence
a constant wave number.
\end{proof}
\end{proposition}
\begin{definition}
A rotation of a multiplicative wave number ${\bf w}(f,g)$ is its product
${\bf w}(f,g+h)={\bf w}(0,h){\bf w}(f,g)$ 
with a constant wave number ${\bf w}(0,h)$.
\end{definition}
\begin{proposition}
The constant wave numbers
form the rotation group of a multiplicative wave number 
${\bf w}(\tfrac{m}{n},\tfrac{p}{q})$.
\end{proposition}
\begin{definition}
A wave number ${\bm \omega}=
{\bf w}(\tfrac{m}{n},\tfrac{p}{q})$ is simple iff $q=n$.
\end{definition}
\begin{proposition}
\label{theorem:simple}
The set of simple multiplicative wave numbers forms
a proper subgroup of the group of multiplicative wave numbers.
\begin{proof}
The set of simple multiplicative wave numbers is closed under the product operation while the identity and the inverse 
of a simple wave number are both simple. 
\end{proof} 
\end{proposition}
The multiplicative wave numbers have a multiplicative norm
\begin{theorem}
\label{theorem:mag1}
The norm
$\|{\bf w}(f,g)\|=1 \ \forall \ {\bf w}(f,g)$ 
\begin{proof}
By Definition \ref{def:norm}
the norm of the periodic sequence ${\bf w}(f,g)$ is
\begin{eqnarray}	
\|  {\bf w}(f,g) \| 
=	\sqrt{ \frac { \underset {n}{\Sigma}
({\bf w}(f,g) \overline{\bf w}(f,g))}   {n}   }
=	\sqrt{ \frac { \underset {n}{\Sigma}1}   {n}   }
=	\sqrt{ \frac {n}   {n}   }
=1.
\end{eqnarray}  
in which $n$ is the period.
\end{proof} 
\end{theorem}
while by Theorem \ref{multiplicative}
the multiplicative property of the magnitude is represented in
\begin{corollary} 
\begin{eqnarray}
\lVert {\bf w}(f,g) {\bf w}(p,q) \rVert =\lVert {\bf w}(f+p,g+q)\| =1
=1.1=\|{\bf w}(f,g)\|\|{\bf w}(p,q) \|
\end{eqnarray}   
\end{corollary} 
\begin{remark}
Interpretations of the multiplicative wave number 
${\bf w}(\tfrac{m}{n},0)$ include one 
in which it is viewed it as a periodic one-dimensional 
string of length $n$ that has
one of two directions.
Under this interpretation the phases of a wave number 
are viewed as co-existing rather than existing
in some sequential sense.
The second interpretation is one in which 
it is viewed as a unit circle with $n$ 
uniformly-spaced points on the circumference, two polarities,
and a rotation rate proportional 
to $\tfrac{m}{n}$. 
The product and sum of two strings under the first interpretation
may then be viewed as the
element-wise sum and product of a concatenation of m strings
of period n and n strings of period m.
Under the second interpretation, 
the sums and products may be viewed  as
occurring when two values coincide at a fixed point of
their two unit circles.
\end{remark}
\section{The Commutative Ring $\mathbb{W}_A(\oplus, \otimes)$ 
of Additive Wave Numbers}
It is natural to extend 
the set ${\bf W}_M$ of multiplicative wave numbers
to include wave numbers that are closed under the addition operator $\oplus$ as well as under the operators  
$\big\{{M}_R, {M}_I, { R}_n,\otimes \big\}$. 
The extended set of wave numbers, denoted by ${\bf W}_A $,
therefore includes a subset of elements that correspond to sums of sequences
representing the $n$th roots of unity.
\subsection{The Definition and Representation of Additive Wave Numbers\\}
\begin{definition} 
\label{def:additive}
A wave number ${\bm \omega}$
is additive and ${\bm \omega}\ \epsilon\ {\bf W}_A $ iff
\begin{eqnarray}
\label{eq:awn}
{\bf W}_A ={\bf W}_U \cup
{\widehat{M}}_R ({\bf W}_A)  \cup 
{\widehat{M}}_I ({\bf W}_A)  \cup  
{\widehat{ R}}_{N} ({\bf W}_A)  \cup  
{\bf W}_A\widehat{\otimes}{\bf W}_A \cup  
{\bf W}_A\widehat{\oplus} {\bf W}_A
\end{eqnarray}
in which ${\bf W}_U$ is the set of unitary wave numbers.
\end{definition}
It follows from Proposition \ref{def:multigen1}
that the additive wave numbers are the closure of ${\bf M}_U$ with respect
to the operators
${\bf S}_M= \big\{{M}_R, {M}_I, { R}_N,\otimes, \oplus \big\}$.
It follows from Definition \ref{def:additive} that the unitary wave numbers 
are additive and from Definition \ref{def:multiplic}
that the multiplicative wave numbers are additive.

In order to find a representation for the additive wave numbers
that generalizes the representation of Theorem \ref{thm:ratform}
for multiplicative wave numbers, it is necessary to find a representation
for the sum of a finite number
of multiplicative wave numbers
\begin{eqnarray}
\label{genform}
S_N({\bf F}_1,...,{\bf F}_N)={\overset {N} {\underset {j=1} 
\oplus }
e^{i {\bf F}_j}   }, \ \  {\bf F}_j =2 \pi (f_j{\bm \xi} \oplus g_j {\bf 1}).
\end{eqnarray}
To that end, it is useful to identify special subsets of 
${\bf W}_A$ that include the constant integer wave numbers ${\bf W}_{INT}=
\Big\{ {\bf ..,-2,-1,0,1,2,...}  \Big\}$.
\begin{proposition} 
\label{def:nwn}
The integer wave numbers 
$  {\bf W}_{INT}
\subset {\bf W}_A $
\begin{proof}
By proposition \ref{prop:i}, 
$ \big\{\pm{\bf 1}\big\} \subset
{\bf W}_M$ and since ${\bf W}_M \subset
{\bf W}_A$ it follows that
$ {\bf1} \oplus (-{\bf 1})=
 {\bf1} \ominus {\bf 1}={\bf 0}
\ \epsilon \   {\bf W}_A $. Hence
\begin{eqnarray}
\label{eq:awn}
{\bf W}_{INT} = \ \big\{  {\bf 1, -1} \big \} \cup  
{\bf W}_{INT}\widehat{\oplus} {\bf W}_{INT}
\end{eqnarray}
\end{proof}
\end{proposition}
defines the integer wave numbers.

It is also useful to show that certain sequences of trigonometric functions are
additive wave numbers.
To this end,
one may decompose the exponential representations of multiplicative
wave numbers into their
real and imaginary components by applying Euler's formula 
in an element-wise manner
\begin{eqnarray}
e^{i 2\pi (f{\bm \xi} \oplus g{\bf 1})} = cos (2\pi (f{\bm \xi}\oplus g{\bf 1}))+{\bf i} sin(2\pi (f{\bm \xi}\oplus g{\bf 1}))
\end{eqnarray}
and defining the sequences
\begin{eqnarray}
{\bf cos}(f,g)  =
\cos(2\pi(f{\bm \xi}\oplus g{\bf 1})); \ \ \ 
{\bf sin}(f,g) =
\sin(2\pi(f{\bm \xi}\oplus g{\bf 1}))
\end{eqnarray}
which,
together with the definitions of addition $\oplus$ and subtraction $\ominus$,
lead to
\begin{lemma}
\label{th:sum} 
The sum and difference
of two multiplicative wave numbers
${\bf w}(f_1,g_1)$ and ${\bf w}(f_2,g_2)$ 
may be represented by the expressions
\begin{eqnarray}
\label{2wavenus}
{\bf w}(f_1,g_1)\oplus {\bf w}(f_2,g_2)
&=&
{\bf 2} {\bf cos} \big(\tfrac{f_1-f_2}{2},\tfrac{g_1-g_2}{2}\big)
 {\bf w}\big(\tfrac{f_1+f_2}{2},\tfrac{g_1+g_2}{2}\big)
 \\
{\bf w}(f_1,g_1)\ominus {\bf w}(f_2,g_2)
&=&
{\bf 2} {\bf i}{\bf sin} \big(\tfrac{f_1-f_2}{2},\tfrac{g_1-g_2}{2}\big)
  {\bf w}\big(\tfrac{f_1+f_2}{2},\tfrac{g_1+g_2}{2}\big) 
\nonumber
\end{eqnarray}
\begin{proof}
Setting ${\bm F}_1=2\pi(f_1{\bm \xi}\oplus g_1{\bf 1})$ 
and ${\bm F}_2= 2\pi(f_2{\bm \xi}\oplus g_2{\bf 1})$,
the sum ${\bf w}(f_1,g_1)\oplus {\bf w}(f_2,g_2)$
of two multiplicative wave numbers becomes
\begin{eqnarray}
\label{eq:proofseq}
e^{i{\bm F}_1}\oplus e^{i{\bf F}_2} 
&=&cos({\bm F}_1)\oplus i sin ({\bm F}_1) \oplus cos ({\bm F}_2) 
\oplus i sin ({\bm F}_2) \\
&=& (cos({\bm F}_1)\oplus cos ({\bm F}_2) )
\oplus  (sin ({\bm F}_1)\oplus  sin ({\bm F}_2) )\nonumber \\
&=& {\bf 2}\cos(\tfrac{{\bm F}_1\oplus {\bm F}_2}{2})\cdot
\cos(\tfrac{{\bm F}_1 \ominus {\bm F}_2}{2})
\oplus i{\bf 2}\sin(\tfrac{{\bm F}_1\oplus {\bm F}_2}{2})\cdot
\cos(\tfrac{{\bm F}_1 \ominus {\bm F}_2}{2})\nonumber \\
&=&{\bf 2}\cos(\tfrac{{\bm F}_1 \ominus {\bm F}_2}{2})\cdot
(\cos(\tfrac{{\bm F}_1\oplus {\bm F}_2}{2})
 \oplus i\sin(\tfrac{{\bm F}_1\oplus {\bm F}_2}{2}))\nonumber \\
&=&{\bf 2}\cos(\tfrac{{\bm F}_1 \ominus {\bm F}_2}{2})\cdot 
e^{\tfrac{i ({\bm F}_1\oplus
{\bm F}_2)}{2}}\nonumber 
\end{eqnarray}
in which ${\bf 2}$ is a constant integer wave number.
The case for  subtraction of wave numbers follows from analogous reasoning.
\end{proof}
\end{lemma}
It follows immediately from Lemma \ref{th:sum} and Theorem
\ref{thm:ratform} that one has
\begin{lemma}
\label{2sin}
${\bf 2}{\bf cos}(f,g)$ and ${\bf 2}{\bf sin}(f,g)$  are additive wave numbers.
\begin{proof}
On post-multiplying each of the wave numbers
in Equation (\ref{2wavenus}) by the conjugate wave number
$ \overline{  {\bf w}}\big(\tfrac{f_1+f_2}{2},\tfrac{g_1+g_2}{2}\big)$,
distributing $\otimes$ over $\oplus$, and multiplying the second equation
by $-{\bf i}$, one obtains
\begin{eqnarray}
{\bf 2} {\bf cos}\big(\tfrac{f_1-f_2}{2}, \tfrac{g_1-g_2}{2}\big)&=&
{\bf w}\big(\tfrac{f_1-f_2}{2},\tfrac{g_1-g_2}{2}\big)\oplus
 {\bf w}\big(\tfrac{f_2-f_1}{2},\tfrac{g_2-g_1}{2}\big)
\\
 {\bf 2}{\bf sin} \big(\tfrac{f_1-f_2}{2}, \tfrac{g_1-g_2}{2}\big)&=&
{\bf w}\big(\tfrac{f_1-f_2}{2},\tfrac{g_1-g_2}{2}\big)\ominus
 {\bf w}\big(\tfrac{f_2-f_1}{2},\tfrac{g_2-g_1}{2}\big) \nonumber
\end{eqnarray}
which, by Definition \ref{def:additive} are additive wave numbers.
\end{proof}
\end{lemma}
One has the immediate
\begin{corollary}
The mirror images, roots, products,  and sums of
${\bf 2}{\bf cos}(f,g)$ and ${\bf 2}{\bf sin}(f,g)$ are additive wave numbers.
\end{corollary}
One may now prove a theorem that generalizes 
Lemma \ref{th:sum} to a sum of $n$
multiplicative wave numbers, and show how 
such a sum  may be
represented as the product of a multiplicative wave number
${\bf w}(f,g)$ and an amplitude function ${\bf A}$ composed of
products of functions of sine and cosine terms.
In particular, the theorem provides a representation of the sum that
is invariant under the order in which the summands are added. 
It is to be noted that the value of the sum is, by associativity
and commutativity, always independent of this order.
\begin{theorem}
\label{thm:wnform11}
The sum $S_N(F_1,...,F_N)={\overset {N} {\underset {j=1} \Sigma }
e^{i F_j}   }$ 
has the representation
\begin{eqnarray}
S_N(F_1,...,F_N)=
{A}_N(F_1,...,F_N) e^{i \big(\overset {N} {\underset {j=1} \Sigma }\frac{ F_j}{N}\big)}
\end{eqnarray}
in which the functions $A_N()$ are defined recursively for all values $N \  \epsilon\ \mathbb{N}$:
\begin{eqnarray}
A_1(F_1)&=&1;\\
A_2(F_1,F_2)&=&2\cos(\tfrac{F_1-F_2)}{2}; \nonumber\\
A_3(F_1,F_2,F_3)&=&
2(A^{\tfrac{1}{2}}_{21} A^{\tfrac{1}{2}}_{22}A^{\tfrac{1}{2}}_{23})^\frac{1}{3}
\Big(2\cos\big(\tfrac{F_2+F_3-2F_1}{4}-iln A^\frac{1}{2}_{21}\big)  \nonumber\\ 
&\ &2\cos\big(\tfrac{F_3+F_1-2F_2}{4}-iln A^{\frac{1}{2}}_{22}\big)
2\cos\big(\tfrac{F_1+F_2-2F_3}{4}-iln A^{\frac{1}{2}}_{23}\big)
\Big)^\frac{1}{3}\nonumber\\
\vdots\nonumber\\
A_{N+1}(F_1,...,F_{N+1})&=&
2\Big({\overset {N+1} {\underset {m=1} \Pi }
A^{\tfrac{1}{2}}_{Nm}}
\Big)^{\frac{1}{N+1}}
\Bigg(
{
\overset {N+1}
{\underset {m=1}
 \Pi }}
{{2\cos\Big({\tfrac{{\overset {N}  {  \underset {j\ne m} {\underset {j=1} 
 \Sigma }}{ F_j}}-{N}F_m}{2N}}}-ilnA^{\tfrac{1}{2}}_{Nm}}\Big)\Bigg)^{\tfrac{1}
{(N+1)}}\nonumber \\
&=&
4
\Bigg(
{
\overset {N+1}
{\underset {m=1}
 \Pi }}A^{\tfrac{1}{2}}_{Nm}
{{\cos\Big({\tfrac{{\overset {N}  {  \underset {j\ne m} {\underset {j=1} 
 \Sigma }}{ F_j}}-{N}F_m}{2N}}}-ilnA^{\tfrac{1}{2}}_{Nm}}\Big)\Bigg)^{\tfrac{1}
{(N+1)}}\nonumber
\end{eqnarray}
in which $A_{Nm}$ indicates a value of $A_N$ 
corresponding to one of $N$ orderings associated
with the $N$ values $F_1,...,F_N$. 
\begin{proof} 
By induction. For the case
$N=1$,
$S(F_1)=e^{i F_1}$.
For $N=2$,  it follows from Lemma \ref{th:sum}
that $S(F_1,F_2) = e^{i{F_1}}+ e^{i{F_2}}=  
2\cos(\tfrac{(F_1-F_2)}{2}e^{i\tfrac{(F_1+F_2)}{2}}$.
The case for $N=3$ provides useful insight into the general case
for $N$ arbitrary by demonstrating how one must
write the sum $S_3$ of three exponential terms 
in three distinct ways, apply Lemma \ref{th:sum} 
in a pairwise manner to the exponential terms,
and then form the one third root of their product
in order to obtain a unique representation of the sum $S_3$. 
One first obtains
\begin{eqnarray}
S_{31}&=&(e^{i{F_2}}+ e^{i{F_3}})+ e^{i{F_1}} =A_{21} e^{i\tfrac{(F_2+F_3)}{2}}+ e^{i{F_1}} 
=e^{i\big(\tfrac{(F_2+F_3)}{2}-ilnA_{21}\big)}+ e^{i{F_1}}
\nonumber\\
&=&
2 cos \big(\tfrac{(F_2+F_3)/2-F_1-ilnA_{21})}{2}\big)e^{i\big(\tfrac{(F_2+F_3)/2+F_1}{2}\big)}A_{21}^{1/2}\\
S_{32}&=&(e^{i{F_3}}+ e^{i{F_1}})+ e^{i{F_2}}  =A_{22} e^{i\tfrac{(F_3+F_1)}{2}}+ e^{i{F_2}}
=e^{i\big(\tfrac{(F_3+F_1)}{2}-ilnA_{22}\big)}+ e^{i{F_2}}\nonumber\\
&=&
2 cos \big(\tfrac{(F_3+F_3)/2-F_2-ilnA_{22})}{2}\big)e^{i\big(\tfrac{(F_3+F_1)/2+F_2}{2}\big)}A_{22}^{1/2}\nonumber\\
S_{33}&=&(e^{i{F_1}}+ e^{i{F_2}})+ e^{i{F_3}} 
=A_{23} e^{i\tfrac{(F_1+F_2)}{2}}+ e^{i{F_3}} 
=e^{i\big(\tfrac{(F_1+F_2)}{2}-ilnA_{23}\big)}+ e^{i{F_3}}\nonumber\\
&=&
2 cos \big(\tfrac{(F_1+F_2)/2-F_3-ilnA_{23})}{2}\big)e^{i\big(\tfrac{(F_1+F_2)/2+F_3}{2}\big)}A_{23}^{1/2}\nonumber
\end{eqnarray}
in which $S_{31}=S_{32}=S_{33}=S$ and $A_{21}= 2 cos(\tfrac{F_2-F_3}{2})$,
$A_{22}= 2 cos(\tfrac{F_3-F_1}{2})$, and 
$A_{23}=2 cos(\tfrac{F_1-F_2}{2})$.
One may then expand the cosine terms into two exponentials which may be simplified
by multiplying each exponential term through 
by the term ($A_{23}^{1/2}, A_{21}^{1/2}, or
A_{22}^{1/2}$) occurring on its RHS, to obtain
\begin{eqnarray}
S_3
&=&
\Big(e^{i\tfrac{(F_1+F_2)/2-F_3)}{2}}  A_{23}+
e^{-i\tfrac{(F_1+F_2)/2-F_3)}{2}}\Big)
e^{i\big(\tfrac{(F_1+F_2)/2+F_3}{2}\big)}
\\
S_3
&=&
\Big(e^{i\tfrac{(F_2+F_3)/2-F_1)}{2}} A_{21}+
e^{-i\tfrac{(F_2+F_3)/2-F_1)}{2} }\Big)
e^{i\big(\tfrac{(F_3+F_1)/2+F_2}{2}\big)}\nonumber \\
S_3
&=&
\Big(e^{i\tfrac{(F_3+F_1)/2-F_2)}{2}} A_{22}+
e^{-i\tfrac{(F_3+F_1)/2-F_2)}{2}} \Big)
e^{i\big(\tfrac{(F_1+F_2)/2+F_3}{2}\big)}\nonumber
\end{eqnarray}
Factoring each equation by its first term,
forming the product of the three equations, 
and taking the $1/3$ root, one obtains
\begin{eqnarray}
\label{eq:factor}
S_3
=&\Big(&
A_{23}+e^{-i\big(\tfrac{(F_1+F_2)-2F_3}{2}\big)}\Big)^\frac{1}{3}
\Big(
A_{21}+e^{-i\big(\tfrac{(F_2+F_3)-2F_1}{2}\big)}\Big)^\frac{1}{3}
\nonumber\\
&\Big(&A_{22}+e^{-i\big(\tfrac{(F_3+F_1)-2F_2}{2}\big)}\Big)^\frac{1}{3} \Big(e^{i\big(\tfrac{F_1+F_2+F_3}{3}\big)}\Big)
\end{eqnarray}
since the product of the factored terms
\begin{eqnarray}
\label{eq:useful}
e^{i\tfrac{(F_1+F_2)/2-F_3)}{2}} e^{i\tfrac{(F_2+F_3)/2-F_1)}{2}} 
e^{i\tfrac{(F_3+F_1)/2-F_3)}{2}} =1.
\end{eqnarray}
One may now apply Lemma \ref{th:sum}, on noting that
$A_{jk}=e^{i (-i ln A_{jk})}$, to each of the first three terms 
in equation (\ref{eq:factor}) to obtain
 \begin{eqnarray}
\label{eq:factor1}
S_3
&=&
4\Big(
 A^\frac{1}{2}_{21}\cos\big(\tfrac{F_2+F_3-2F_1}{4}-iln A^\frac{1}{2}_{21}\big) 
\ A^{\frac{1}{2}}_{22}\cos\big(\tfrac{F_3+F_1-2F_2}{4}-iln A^{\frac{1}{2}}_{22}\big)\nonumber \\
&\ & A^{\frac{1}{2}}_{23}\cos\big(\tfrac{F_1+F_2-2F_3}{4}-iln A^{\frac{1}{2}}_{23}\big)
\Big)^\frac{1}{3} \ e^{i\big(\tfrac{F_1+F_2+F_3}{3}\big)} 
\end{eqnarray}
following a second application of equation (\ref{eq:useful}).
For the general case $N+1$, one may
assume the inductive hypothesis
\begin{eqnarray}S_N=
{\overset {N} {\underset {j=1} \Sigma }e^{i F_j}}=
{
A_{N} e^{i \big(\overset {N} {\underset {k=1} \Sigma }\frac{ F_j}{N}\big)}},
\hspace{.3in} 
\end{eqnarray}
for arbitrary $N$ and follow the steps in the proof for the case $N=3$.
One first
expresses the sum of $N+1$ exponential terms in $N+1$ distinct ways as
\begin{eqnarray}
S_{N+1,m}&=&
{\overset {N+1} {\underset {j=1} \Sigma }e^{i F_j}}
={\overset {N}  {  \underset {j\ne m} {\underset {j=1} 
 \Sigma }}e^{i F_j}}+e^{iF_m}, \hspace{.3in}m=1,N+1
\end{eqnarray}
in which a given value of $m$ represents one of the alternative representations
and applies the inductive hypothesis to each representation, obtaining
\begin{eqnarray}
S_{N+1,m}
=
A_{N,m}e^{i \big(\overset {N+1}  {  \underset {j\ne m} {\underset {j=1} 
 \Sigma }}\frac{ F_j}{N}\big)}
+e^{iF_m}
=
e^{i \big(\overset {N+1}  {  \underset {j\ne m} {\underset {j=1} 
 \Sigma }}\frac{ F_j}{N}-ilnA_{N,m}\big)}
+e^{iF_m}
, \hspace{.3in}  m = 1,N+1
\end{eqnarray}
noting that, in general, $A_{N+1,m}\ne A_{N+1,n}$ for $m\ne n$.
It follows from Lemma \ref{th:sum} for the sum of two exponential terms
that
\begin{eqnarray}
S_{N+1,m}
=
2 \cos\Bigg(\Big({{\overset {N+1}  {  \underset {j\ne m} {\underset {j=1} 
 \Sigma }}\tfrac{ F_j}{N}}-F_m-iln A_{N,m}}\Big)\Big/{2}\Bigg)
e^{i\Bigg(\Big({{\overset {N+1}  {  \underset {j\ne m} {\underset {j=1} 
 \Sigma }}\tfrac{ F_j}{N}}+F_m-iln A_{N,m}}\Big)\Big/{2}\Bigg)}
\end{eqnarray}
and, as in the case $N=3$, one may transform the cosine term into the sum of
two exponential terms and simplify the resulting equation to obtain
\begin{eqnarray}
S_{N+1,m}
&=&\Bigg(A_{N,m}e^{i\Big({{\overset {N+1}  {  \underset {j\ne m} {\underset {j=1} 
 \Sigma }}{ F_j}-{N}}F_m}\Big)\Big/{2N}}+
e^{-i\Big({{\overset {N+1}  {  \underset {j\ne m} {\underset {j=1} 
 \Sigma }}{ F_j}-{N}}F_m}\Big)\Big/{2N}}\Bigg)
e^{i\Big({{\overset {N}  {  \underset {j\ne m} {\underset {j=1} 
 \Sigma }}{ F_j}}+{N}F_m}\Big)\Big/{2N}} \nonumber \\
&= &
e^{i\Big({{\overset {N}  {  \underset {j\ne m} {\underset {j=1} 
 \Sigma }}{ F_j}}-{N}F_m}\Big)\Big/{2N}}
\Bigg(A_{N,m}+
e^{-i\Big({{\overset {N+1}  {  \underset {j\ne m} {\underset {j=1} 
 \Sigma }}{ F_j}-{N}}F_m}\Big)\Big/{N}}\Bigg)
e^{i\Big({{\overset {N}  {  \underset {j\ne m} {\underset {j=1} 
 \Sigma }}{ F_j}}+{N}F_m}\Big)\Big/{2N}} 
\end{eqnarray}
After multiplying together the $N+1$ representations of $S_{N+1}$
and taking the $(N+1)th$ root of the product, one obtains
\begin{eqnarray}
S_{N+1}
={\overset {N+1} {\underset {m=1} \Pi }}
\Bigg({{e^{i(-ilnA_{Nm})}+e^{-i\Big({{\overset {N}  {  \underset {j\ne m} {\underset {j=1} 
 \Sigma }}{ F_j}}-{N}F_m}\Big)\Big/{N}}}}\Bigg)^{\tfrac{1}{(N+1)}}
 e^{i\Big(\frac{\overset {N+1} {   {\underset {j=1} 
\Sigma }} F_j}{N+1}\Big)}
 \hspace{.3in}m=1,N+1
\end{eqnarray}
Applying Theorem \ref{th:sum} to the two terms in each of the product
terms
\begin{eqnarray}
S_{N+1}
&=&2\Big({\overset {N+1} {\underset {m=1} \Pi }
A^{\tfrac{1}{2}}_{Nm}}
\Big)^{\frac{1}{N+1}}
\Bigg(
{
\overset {N+1}
{\underset {m=1}
 \Pi }}
{{2\cos\Big({\tfrac{{\overset {N}  {  \underset {j\ne m} {\underset {j=1} 
 \Sigma }}{ F_j}}-{N}F_m}{2N}}}-ilnA^{\tfrac{1}{2}}_{Nm}}\Big)\Bigg)^{\tfrac{1}
{(N+1)}}
 e^{i\Big(\tfrac{\overset {N+1} {   {\underset {j=1} 
\Sigma }} F_j}{N+1}\Big)} 
\nonumber \\
&=&
4
\Bigg(
{
\overset {N+1}
{\underset {m=1}
 \Pi }}A^{\tfrac{1}{2}}_{Nm}
{{\cos\Big({\tfrac{{\overset {N}  {  \underset {j\ne m} {\underset {j=1} 
 \Sigma }}{ F_j}}-{N}F_m}{2N}}}-ilnA^{\tfrac{1}{2}}_{Nm}}\Big)\Bigg)^{\tfrac{1}
{(N+1)}} e^{i\Big(\tfrac{\overset {N+1} {   {\underset {j=1} 
\Sigma }} F_j}{N+1}\Big)}\nonumber \\
&=&A_{N+1} 
 e^{i\Big(\tfrac{\overset {N+1} {   {\underset {j=1} 
\Sigma }} F_j}{N+1}\Big)}
\end{eqnarray} \qed
\end{proof}
\end{theorem}
One notes that Theorem \ref{thm:wnform11} 
applies to sums of exponential functions
for which the argument $F_j, j=1,N$ may have different forms,
including that of a real or complex scalar or
that of a sequence of scalars. In the latter case, one notes that a 
sum of exponentials whose arguments are sequences
is defined in a component-wise manner and hence the theorem applies.
Since the latter case includes multiplicative wave numbers, 
which are defined by $ {\bf F}_j =2 \pi (f_j{\bm \xi} \oplus g_j {\bf 1})$
one immediately has
\begin{corollary}
\label{cor:genform1}
A finite sum of multiplicative wave numbers 
may be represented as
a product ${\bf A}{\bf w}$ 
of an amplitude function $\bf A$
that is representable as a product
of sine and cosine factors and a multiplicative wave number 
${\bf w}$.
\begin{proof}
On assuming
$ {\bf F}_j =2 \pi (f_j{\bm \xi} \oplus g_j {\bf 1})$ in Theorem \ref{thm:wnform11} 
it follows that
\begin{eqnarray}
\overset {N} {\underset {j=1} \oplus } e^{i {\bf F}_j}    
&=&
{A}_N({\bf F}_1,...,{\bf F}_N) e^{i \big(\overset {N} {\underset {j=1} \oplus }\frac{{\bf F}_j}{N}\big)} \nonumber\\
&=&{A}_N({\bf F}_1,...,{\bf F}_N) 
{\bf w} \Big(   \tfrac   { \overset {N} {\underset {j=1} \Sigma} f_j  } {N},
 \tfrac   { \overset {N} {\underset {j=1} \Sigma }  g_j} {N} \Big) 
\end{eqnarray}
\end{proof}
\end{corollary}
Before generalizing the polar representation of Theorem \ref{thm:ratform}
from multiplicative wave numbers to additive wave numbers,
it is useful to make
\begin{definition}
The amplitude ${\bf A}$
of the wavenumber ${\bf A}  {\bf w} \ \epsilon\ {\bf W}_A$
is in w-free form iff every factor ${\bf a}$ of ${\bf A}$
such that ${\bf a}\ \epsilon \  {\bf W}_M$ is
combined into ${\bf w}$.
\end{definition}
\begin{theorem}
\label{thm:genform1}
An additive wave number ${\boldsymbol \omega}$ has a general
representation
as the product of an w-free amplitude function 
${\bf A}$ that is a product of sine and cosine functions and their mirror images and roots and a multiplicative wave number
$
{\bf w}(f,g)$
\begin{proof}
One first notes that multiplicative wave numbers,
together with their reflection 
numbers for which the amplitude function $A=1$.
As shown in Corollary \ref{cor:genform1} 
a finite sum of multiplicative wave numbers has the general form
${\bf A}{\bf w}(f,g) $
This form is invariant under applications of the reflection, root,
product, and sum operations ${M}_R, {M}_I, {\it R_N}, \otimes, \oplus$:
\begin {enumerate}
\item
for the mirror image operators 
\begin{eqnarray}
{M}_R({\bf A}{\bf w}(f,g)) &=&\bar{\bf A}\bar{\bf w}(f,g)\\
{M}_I ({\bf A}{\bf w}(f,s)) &=&
-\bar{\bf A}\bar{\bf w}(f,g)\nonumber 
\end{eqnarray}
\item
for the product operator
\begin{eqnarray}
{\bf A}_1{\bf w}(f_1,g_1){\bf A}_2{\bf w}(f_2,g_2) &=&
({\bf A}_1{\bf A}_2 ){\bf w}(f_1,g_1){\bf w}(f_2,g_2)  \nonumber \\
&=&
{\bf A}_{12} {\bf w}(f_1+f_2,g_1+g_2) 
\end{eqnarray}
\item
for the root operator
\begin{eqnarray}
{\it R}_n({\bf A}{\bf w}(f,g)) &=&{\it R}_n({\bf A}) 
 {\it R}_n( {\bf w}(f,g)) 
\end{eqnarray}
\item
for the sum operator, the general form of the representation is proved in Corollary \ref{cor:genform1}.
\end{enumerate}
Hence there is no sequence of the operations 
${M}_R, {M}_I, {\it R}_n, \otimes,  \oplus$ 
leading to an additive wave number
that cannot be represented in the form 
${\bf A}{\bf w}(f,g) $.
Furthermore
all amplitudes ${\bf A}$ may be assumed to be in w-free form
by the
commutativity of $\otimes$, the invariance
of multiplicative wave numbers under the mirror image, root,
and product operations, and the finiteness of the
periodicity of wave numbers.
\qed
\end{proof}
\end{theorem}
\subsection{Some Elementary Properties of the Additive Wave Numbers\\}
The elementary properties of additive wave numbers 
that are considered here 
include
their periodicity, their ring and module structure, and their multiplicative norm.
The periodicity of an additive wave number follows from
Theorem \ref{theor:period}.
\begin{theorem}
The set $\bf{W}_	A$ of additive wave numbers forms 
a commutative ring $\mathbb{W}(\oplus,\otimes)$.
\begin{proof}
It follows from the definition of an additive wave number
that the set ${\bf W}_A$ is closed under $\oplus$.
The wave number
\begin{eqnarray}
{\bf w}(f,g) \ominus{\bf w}(f,g) =\Big\{ 2i sin(0) e^{2\pi i (f\xi+s)}\Big\}
=\Big\{0\Big\}= {\bm 0}
\end{eqnarray}
is both an additive wave number and
the identity of the sum operator
$\oplus$
since
\begin{eqnarray}
A {\bf w}(f,g) \oplus {\bm 0}=
\Big\{ a(\xi)  e^{2\pi i( f\xi+s )}+0  
 \Big\}
=\Big\{ a(\xi)    e^{2\pi i( f\xi+s )}  
 \Big\}
=A {\bf w}(f,g) .
\end{eqnarray}
It follows from the definition of the additive wave numbers
that if $A {\bf w}(f,g) $ exists then the additive wave number $-A {\bf w}(f,g) )$ also  exists. The relation
\begin{eqnarray}
A {\bf w}(f,g) \oplus (-A {\bf w}(f,g)  )=A {\bf w}(f,g) 
\ominus A {\bf w}(f,g) 
= {\bm 0}.
\end{eqnarray}
follows from the definition of the equality of countable sequences and hence
every discrete wave number has a unique additive inverse.
Finally it follows from the 
commutativity and associativity of $\oplus$ that 
$ {\bf W}_A(\oplus,\otimes)$is an Abelian additive group.
Since, however, the product operator $\otimes$ does not have
a multiplicative inverse for any amplitude apart from ${\bf 1}$,
the set $ {\bf W}_A$ forms a commutative ring.
\qed
\end{proof}
\end{theorem}
One may show further that the set ${\bf W}_A$ has a module-like structure 
in which the scalars of a ring are replaced by a set ${\bf F}_A$ of amplitudes.
\begin{proposition}
\label{prop11}
${\bf A} {\bf w}\ \epsilon \ {\bf W}_A \Rightarrow {\bf A}\ \epsilon \ {\bf W}_A$
\begin{proof}
${\bf w}\ \epsilon \ {\bf W}_M \Rightarrow {\bf w}\ \epsilon \ {\bf W}_A
 \Rightarrow \overline{{\bf w}}\ \epsilon \ {\bf W}_A. \ Hence \ 
{\bf A} {\bf w}\overline{\bf w}={\bf A}\ \epsilon\ {\bf W}_A.$
\end{proof}
\end{proposition} 
In showing that ${\bf W}_A$ has a module-like structure, it is useful to make
\begin{definition}
${\bf F}_A=\big\{ {\bf A} \mid \forall {\bf w}\ \epsilon\ {\bf W}_M, 
{\bf A} {\bf w} \ \epsilon\ {\bf W}_A \big\}$.
\end{definition}
\begin{proposition}
\label{th}
${\bf F}_A$ is a proper subset, but not a subring, of $ {\bf  W}_A$
that is closed under $\otimes$ and has ${\bf F}_A\cap {\bf W}_M ={\bf 1}$.
\begin{proof}
By Proposition \ref{prop11},
${\bf F}_A$ is a subset of ${\bf W}_A$.
For ${\bf w}(f,g)
\ \epsilon\  {\bf W}_M,\\
{\bf w}(f,g) \notin {\bf F}_A $
by the assumption that ${\bf A}$ is ${\bf w}$-free. 
Hence ${\bf F}_A$ is a proper subset of ${\bf W}_A $.
Since ${\bf 1} \ \epsilon\ {\bf F}_A$ by definition, 
${\bf F}_A\cap {\bf W}_M ={\bf 1}$.
For ${\bf A}= cos(a), {\bf A}' =i sin(a)$, one has ${\bf A} \oplus {\bf A}' =e^{ia}$ 
which is not permitted. Hence
${\bf F}_A$ is not an additive group.
For ${\bf A}_{1}, {\bf A}_{2}\ \epsilon \ {\bf F}_A$,
${\bf A}_{1} {\bf w}_2\ \epsilon \ {\bf W}_A$, 
and hence ${\bf A}_1, {\bf A}_{2} {\bf w}_1, {\bf A}_{2}{\bf w}_1\ \epsilon \ {\bf W}_A$
and by the closure of ${\bf W}_A$ under $\otimes$, 
${\bf A}_{1} {\bf A}_{2}\ \epsilon \ {\bf F}_A$
\end{proof}
\end{proposition}
Hence one has
\begin{proposition}
\label{th1}
${\bf W}_A ={\bf {\bf F}}_A \hat{\otimes}  {\bf W}_M$
in which $ \hat{\otimes}$ is the set product operator.
\end{proposition}
and
\begin{theorem}
\label{th1}
The multiplicative group ${\bf W}_{M0} ={\bf W}_{M} \cup {\bf 0} $ is a vector space over the set ${\bf F}_A$.
\begin{proof}
The vector ${\bf 0}\ \epsilon\ {\bf W}_{M0}$ by definition.
For all vectors ${\bf w}\ \epsilon {\bf W}_{M0} $ and scalars
${\bf A}\ \epsilon {\bf F}_{A} $, it follows from the properties
of ${\bf W}_A$ that vector addition and
scalar multiplication of vectors are well-defined and satisfy the commutative,
associative, and distributive axioms of a vector space.
\end{proof}
\end{theorem}
The additive wave numbers possess a norm that
is multiplicative.
\begin{theorem}
\label{thm:normadd}
 The norm of the wave number $A{\bf w}$ is
$\|A{\bf w}\|=\sqrt{\frac{ \underset {\lambda}{\Sigma}		
a^2({ \xi})}{\lambda}	}	
$ in which $a({ \xi})$ is a component of the amplitude function 
${\bf A({\bm \xi})}$	
and $\lambda$ is the period of the wave number.	
\begin{proof}
$
 \|A{\bf w}\|=\sqrt{\frac{\underset {\lambda}
{\Sigma} {A {\bf w}} \overline{A {\bf w}}}{\lambda} }
=\sqrt{ \frac{\underset {\lambda}{\Sigma}({\bf A}\bar{\bf A}
{\bf w}\overline{\bf w})}{\lambda}		}
=\sqrt{ \frac{\underset {\lambda}{\Sigma}({\bf A}\bar{\bf A}
)}{\lambda}		}
=\sqrt{\frac{ \underset {\lambda}{\Sigma}		
a^2({ \xi})}{\lambda}	}$					
\end{proof}
\end{theorem}
By Theorem \ref{multiplicative} the norm for periodic sequences is multiplicative,
hence
 \begin{corollary}
The magnitude function for additive wave numbers 
is multiplicative with
\begin{eqnarray}	
\|{\bf A}_1{\bf w}_1\otimes{\bf A}_2{\bf w}_2\|
=\|{\bf A}_1{\bf w}_1\|\|{\bf A}_2{\bf w}_2\|.
\end{eqnarray}	  
\end{corollary}
\begin{example}
The norm of the constant integer wave number ${\bf z}$ is
$\|{\bf z}\|=|{z}|$ 
\end{example}
\begin{proposition}
The set of multiplicative wave numbers $ {\bf W}_M$ 
is not closed
under 
$\oplus$ and
is therefore not an additive group under $\oplus$.
\begin{proof}
The additive wave number 
${\bf 1}\oplus {\bf 1}={\bf 2}$ is not a multiplicative wave number 
since by Theorem \ref{theorem:mag1} any multiplicative wave number
has norm $1$. 
\end{proof} 
\end{proposition}
\subsection{Bases for the Additive Wave Numbers of Period $n$}
\label{sec:bases}
One may construct an orthogonal basis 
for additive wave numbers for each period $n$,
as shown in
\begin{theorem} 
\label{thm:basis}
For every $n\ \epsilon\ 
\mathbb{N}$, 
there exists a set of $n$ orthogonal additive
wave numbers 
of period $n$ whose principal sequences have the form
\begin{eqnarray}
{_1}{\bm u}(\tfrac{1}{n},j)=
\big(0,...,0,n, 0,...,0    \big), \hspace {.2in} j \ \epsilon\  \mathbb{N}
\end{eqnarray}
in which the zeros represent the origin of the complex plane
and
the non-zero term $n$ occurs in the $j$th position.
\begin{proof}
Translating the wave number $ {\bf w}(\tfrac{1}{n},0)$ with each of the
constant wave numbers ${\bf w}(0,\tfrac{j}{n})$ for $j=1,n$ 
one obtains the set of $n$ wave numbers
\begin{eqnarray}
{\bm w_{nj}}={\bf w}(\tfrac{1}{n},0)\ominus{\bf w}(0,\tfrac{j}{n})
\hspace{.3in} j=1,n
\end{eqnarray}
for which the $\xi$th element of the principal sequence
takes the form
\begin{eqnarray}
\label{eq:form1}
{_1}{\bm w_{nj}}(\xi)
&=& \big(e^{\substack{2 \pi i\\ }\frac{\xi}{n}}- e^{\substack{2 \pi i\\ }\frac{j}{n}}\big), \ j=1,n
\end{eqnarray}
in which the diagonal element ${\bf w}_{nj}(j)=0$.
Geometrically, each of these wave numbers represents 
the translation of an n-gon 
by the negative of the vector representing one of its $n$ vertices,
which itself is translated to the origin of the complex plane
in one of the translated n-gons.

One may apply lemma \ref{th:sum} to 
each of the $n^2$ difference terms
\begin{eqnarray}
\label{eq:form1}
{_1}{\bm w_{nj}}(\xi)&=& \big(e^{\substack{2 \pi i\\ }\frac{\xi}{n}}- e^{\substack{2 \pi i\\ }\frac{j}{n}}\big)
=
\Big(2i \sin( \tfrac{2 \pi(\xi-j)}{2n}) e^{\big(\tfrac{2 \pi i (\xi+j)}{2n}\big)}
\Big)
\end{eqnarray}
and then form
the set of $n$ products over $j$
of each distinct combination
of $n-1$ of the $n$ translated wave numbers
(the product of all $n$ results in the wave number ${\bm 0}$).
Since these products define wave numbers
and since the conjugate of a wave number is a wave number, they may
be represented as the conjugate of a wave number 
$\big\{{\bm U}_{n\xi}, \xi \ \epsilon\  \mathbb{Z}\big\}$
whose principal sequence is
\begin{eqnarray}
\label{eq:theset}
{_1}\overline{{\bm U}}_{n\xi} &=&
{_1} {\underset {j\ne \xi}\otimes} {\bf w}_{nj} \\
&=&\Big(0,...,0,
t(\xi),
0,...,0\Big),\hspace{.2in} \xi \ \epsilon\  \mathbb{Z}
\end{eqnarray}
and in which 
\begin{eqnarray}
t(\xi)
&=&
\overset {} {\underset {j\ne k}\Pi}
\Big(2i\sin( \tfrac{2 \pi(\xi-j)}{2n}) e^{2 \pi i \tfrac{(\xi+j)}{2n}}\Big)
 \hspace{.3in} \xi=1,n-1\nonumber \\
&=&
(2i)^{n-1}\Big((-1)^{n-\xi}
\overset {} {\underset {j\ne \xi}\Pi}
\sin( \tfrac{2 \pi(|\xi-j|)}{2n} \Big)
\Big(\overset {} {\underset {j\ne \xi}\Pi}
 e^{2 \pi i  \tfrac{\xi}{2n}}\Big)
\Big(\overset {} {\underset {j\ne \xi}\Pi}
 e^{2 \pi i  \tfrac{j}{2n}}\Big)
 \hspace{.3in} \nonumber \\
&=&
2^{n-1}i^{n-1}\Big((-1)^{n-\xi}
\overset {n-1} {\underset {m=1}\Pi}
 \sin( \tfrac{ \pi m}{n}) \Big)
\Big(
 e^{2 \pi i  \tfrac{\xi(n-1)}{2n}}\Big)
\Big(
 e^{2 \pi i  \tfrac{(n+1)n/2-\xi}{2n}}\Big)
  \nonumber\\
&=&
2^{n-1} e^{2 \pi i \frac{(n-1)n}{4n}}\Big(e^{2 \pi i\frac{2(n-\xi)n}{4n}}
\frac{n}{2^{n-1} }\Big)
\Big(
 e^{2 \pi i \tfrac{(2\xi(n-1) +(n+1)n-2\xi)}{4n}}\Big)
  \nonumber\\
&=&
{n}
 e^{2 \pi i  \tfrac{(n-1)n+2(n-\xi)n+2\xi(n-1) +(n+1)n-2\xi))}{4n}}
\nonumber\\
&=&n
e^{\tfrac{2 \pi i (n^2-\xi)}{n}} \nonumber \\
&=&ne^{-2\pi i \tfrac{\xi}{n}}.
\end{eqnarray}
The term $\overset {n-1} {\underset {m=1}\Pi}
 \sin( \tfrac{ \pi m}{n})$ in the third line follows from the symmetry
of the sine function over the domain $0 \le \theta \le \pi$
and
the fourth line follows from the trigonometric identity
\begin{eqnarray}
\overset {n-1} {\underset {\xi=1}
 \Pi }\sin \Big(\frac{\pi \xi}{n}\Big)=\frac{n}{2^{n-1}}.
\end{eqnarray}
The principal sequences of the
wave numbers may therefore be written as
\begin{eqnarray}
{_1} \overline{{\bm U}}_{n\xi}
&=&n\Big(0,...,0,e^{-2 \pi i \tfrac{\xi
}{n}} ,0,...,0\Big), \hspace{.2in} \xi=1,n
\end{eqnarray}
and one may take its product with the multiplicative wave number
${\bf w}(\tfrac{1}{n},0)=e^{i 2 \pi \tfrac{\bm \xi}{n}}$ to obtain
the additive wave numbers with principal sequence
\begin{eqnarray}
\label{eq:basis}
{_1} {\bm u}(\tfrac{1}{n},j)
&=&\Big(0,...,0,n ,0,...,0\Big) \hspace{.1in}for \  j=1,n
\end{eqnarray} 
in which $n$ occurs in the jth phase. They
may be termed the $n$ phases of ${\bf n}$.
\qed 
\end{proof}
\end{theorem}
\begin{proposition}
The $n$ wave numbers ${\bm u}(\tfrac{1}{n},j),j=1,n$ are mutually
orthogonal under $\otimes$.
\end{proposition}
\section{The Field ${\mathbb W}_I $ of Invertible Wave Numbers}
Since the multiplicative
inverse of an additive wave number ${\bf A }{\bf w}(f,g)$ 
is defined iff ${\bf A }={\bf 1 },$
the additive wave numbers ${\bf W}_A$
form neither a multiplicative group nor a field.
This suggests that one augment the set of operators
that generate the
additive wave numbers with
the inverse operator
${\it I}({\bf A}{ \bf {w}})={\it I}({\bf A}){\it I}({ \bf {w}})
=\frac{  \overline{ \bf {w}} }{{\bf A}}$. 
One notes that ${\it I}$
prohibits as undefined the inverse of any wave number
with one or more zeroes among its elements, generalizing
the prohibition of division by zero in the field of rational numbers.
\subsection{The Definition and Representation of Invertible Wave Numbers\\}
\begin{definition} 
\label{def:invertible}
A wave number ${\bm \omega}$
is invertible and ${\bm \omega}\ \epsilon\ {\bf W}_I $ iff
\begin{eqnarray}
\label{eq:iwn}
{\bf W}_I ={\bf W}_U &\cup&
{\widehat{{M}}_R}({\bf W}_I)  \cup 
{\widehat{{M}}_I}({\bf W}_I)  \cup  {\widehat{{ R}}_n}({\bf W}_I)  \cup  
{\bf W}_I\widehat{\otimes}{\bf W}_I
  \cup
{\bf W}_I\widehat{\oplus}{\bf W}_I \nonumber\\
&\cup&
\widehat{I}({\bf W}_I)
\end{eqnarray}
in which ${\bf W}_U$ is the set of unitary wave numbers.
\end{definition}
It follows from Proposition \ref{def:multigen1}
that the invertible wave numbers are the closure of ${\bf W}_U$ with respect to the operators
${\bf S}_M= \big\{{M}_R, {M}_I, { R}_n,\otimes, \oplus, I \big\}$
while the unitary, multiplicative, and additive wave numbers are invertible.

The inverse operator $I$ allows one to introduce
the quotient operator
\begin{definition}
The
quotient of invertible wave numbers 
${\bm \omega}_1, {\bm \omega}_2 $ is 
\begin{eqnarray}
{\bm \omega}_1 \odiv {\bm \omega}_2 
={\bm \omega}_1\otimes {\it I}( {\bm \omega}_2 ) 
=\frac{{\bm \omega}_1 }{{\bm \omega}_2 }.
\end{eqnarray}
\end{definition}
\begin {example}
The application of the inverse operator to the set of
integer wave numbers ${\bf W}_{Z}$ leads to the set of inverses  
${\bf W}_{Z}^{-1}=\widehat{I}({\bf W}_{Z})$
from which one may define the set of
constant rational wave numbers
\begin{eqnarray}
{\bf W}_Q =
{\bf W}_Z\widehat{\otimes}{\bf W}_Z^{-1}
={\bf W}_Z\widehat{\odiv}{\bf W}_Z
\end{eqnarray}
in which division by ${\bf 0}$ is not permitted.
One may show that they form a field $\mathbb{W}_Q$
having a bijective morphism with $\mathbb Q$.
There is a one-to-one relation between 
the rational invertible wave numbers ${\bf W}_Q$ 
and circles of rational diameter in the complex plane,
with negative rational numbers corresponding to
circles having clockwise polarity.
\end{example}
\begin {proposition}
\label{prop:sincos}
Sine and cosine terms and their inverses are invertible wave numbers.
\begin{proof}
By Theorem \ref{2sin}, ${\bf 2} {\bf \sin ({a})}$ and ${\bf 2} {\bf \cos ({a})}$ 
are additive wave numbers. By $\otimes$-closure, 
${\bf \sin ({a})}={\bf {\tfrac{1}{2}} }({\bf 2} {\bf \sin ({a})})$ and
${\bf \cos({a})}={\bf {\tfrac{1}{2}} }({\bf 2} {\bf \cos ({a})})$
are invertible wave numbers. 
By Definition \ref{def:invertible}, their inverses are invertible wave
numbers.
\end{proof}
\end{proposition}
The invertible wave numbers are well-defined and,
in some sense, "natural" units, as suggested by a
general representation of their sums and differences
that follows from Corollary \ref{nice}.
\begin{theorem}
\label{SandD}
Sums and differences of invertible wave numbers
${\bm \omega}_1=A_1 {\bf w}(f_1,g_1)$, 
${\bm \omega}_2=A_2 {\bf w}(f_2,g_2)$ may be
represented in terms of products and quotients as
\begin{eqnarray}
\label{eq:recurs}
{\bm \omega}_1 \oplus {\bm \omega}_2&=& 2 \cos\Big( \frac{1}{i}
ln\Big(\frac{{\bm \omega}_1}{{\bm \omega}_2}\Big)^{1/2}\Big)
\Big({\bm \omega}_1 {\bm \omega}_2 \Big)^{1/2}\\
{\bm \omega}_1 \ominus {\bm \omega}_2&=& 2i \sin\Big( \frac{1}{i}
ln\Big(\frac{{\bm \omega}_1}{{\bf w}_2}\Big)^{1/2}\Big)
\Big({\bm \omega}_1{\bm \omega}_2 \Big)^{1/2}\nonumber
\end{eqnarray}
\begin{proof} For a sum, one has
\begin{eqnarray}
{\bm \omega}_1 \oplus {\bm \omega}_2&=& {\bf A}_1 {\bf w}(f_1,g_1) \oplus
A_2 {\bf w}(f_2,g_2) \nonumber\\
&=& e^{i{( -iln {\bf A}_1 \oplus 2\pi (f_1{\bm \xi} + g_1{\bf 1}) )}}
\oplus  e^{i{( -iln {\bf A}_2 \oplus 2 \pi(f_2{\bm \xi} + g_2{\bf 1})) }}\nonumber\\
\label{thiseqn}
&=& 2cos(  \tfrac{( -iln {\bf A}_1 \oplus 2\pi (f_1{\bm \xi} + g_1{\bf 1}) )
\ominus( -iln {\bf A}_2 \oplus 2 \pi(f_2{\bm \xi} + g_2{\bf 1}))
}{2}   )\nonumber\\
&\ & \ \   \otimes e^{i     \frac{( -iln {\bf A}_1 \oplus 2\pi (f_1{\bm \xi} + g_1{\bf 1}) )
\oplus( -iln {\bf A}_2 \oplus 2 \pi(f_2{\bm \xi} + g_2{\bf 1}))
}{2}   )   }\nonumber\nonumber\\
&=& 2 cos\Big(\frac{1}{i}      
 ln \Big(   \frac{ {\bf A}_1e^{  2\pi i(f_1{\bm \xi} + g_1{\bf 1}) ) }  }
{ {\bf A}_2  e^{  2\pi i(f_2{\bm \xi} + g_2{\bf 1}) ) }}                                  
   \Big)^{1/2} \Big)\nonumber\\
&\ & \ \   \otimes e^{ln( {\bf A}_1e^{  2\pi i(f_1{\bm \xi} + g_1{\bf 1}) ) }  
 {\bf A}_2  e^{  2\pi i(f_2{\bm \xi} + g_2{\bf 1}) ) }   ) ^{1/2}    }\nonumber\\
&=& 2 \cos\Big( \frac{1}{i}
ln\Big(\frac{{\bm \omega}_1}{{\bm \omega}_2}\Big)^{1/2}\Big)
\Big({\bm \omega}_1 {\bm \omega}_2 \Big)^{1/2}.
\end{eqnarray}
in which the equality
 $2 \pi (f {\bm \xi} + g {\bf 1})=-i ln({\bf w}(f,g)$
 for the multiplicative wave number
${\bf w}(f,g) =e^{i2 \pi (f {\bm \xi} + g {\bf 1})}$
has been used.
An analogous argument applies to the difference.\qed
\end{proof}
\end{theorem}
\begin{corollary}
Occurrences of the $\oplus$ operator in arithmetic expressions 
of invertible wave numbers 
may  be replaced by expressions involving products
and quotients of invertible wave numbers.
\end{corollary}
\begin{corollary}
The sine and cosine functions may be defined in terms of logarithmic 
functions of invertible wave numbers.
\begin{proof}
Dividing the first equation of Theorem \ref{SandD} by the term 
$\Big({\bm \omega}_1 {\bm \omega}_2 \Big)^{1/2}$
and defining ${\bm \rho }= iln({\bf A}w)^{1/2}$ one obtains
$
\label{eq:recurs}
\cos({\bm \rho })= \frac{e^{i{\bm \rho }} \oplus e^{-i{\bm \rho }} }{2}
$.
\end{proof}
\end{corollary}
with the analogous result $
\label{eq:recurs}
\sin({\bm \rho })= \frac{e^{i{\bm \rho }} \ominus e^{-i{\bm \rho }} }{2}
$ for the sine function.

The general form of an invertible wave numbers may now be specified in
\begin{theorem}
\label{thm:genform1}
An invertible wave number ${\bm \omega}$ is representable
as a product ${\bm \omega}={\bf A}{\bf w}(f,g)$ of a multiplicative wave number ${\bf w}(f,g)$ and an invertible amplitude function ${\bf A}$
that is a product of the reflections, powers, products,
and inverses of the sine and cosine functions.
\begin{proof}
By Theorem \ref{SandD},
the sum of invertible wave numbers 
${\bm \omega}_1=A_1 {\bf w}(f_1,g_1)$ and 
${\bm \omega}_2=A_2 {\bf w}(f_2,g_2)$
has the form
\begin{eqnarray}
\label{eq:recurs}
{\bm \omega}_1 \oplus {\bm \omega}_2
&=& 
2 \cos\Big( \frac{1}{i}
ln\Big(\frac{{\bm \omega}_1}{{\bm \omega}_2}\Big)^{1/2}\Big)
\Big({\bm \omega}_1 {\bm \omega}_2 \Big)^{1/2}\\
&=&
\Big(2 \cos\Big( \frac{1}{i}
ln\Big(\frac{{\bm \omega}_1}{{\bm \omega}_2}\Big)^{1/2}\Big)  
\Big( {\bf A}_1 {\bf A}_2\Big)^{1/2} \Big)
 \Big({\bf w}(f_1,g_1) {\bf w}(f_2,g_2) \Big)^{1/2}\nonumber\\
&=&
{\bf A}{\bf w}(f,g)
\end{eqnarray}
with a sine function in place of the cosine in the case of the difference.
The amplitudes ${\bf A}_1, {\bf A}_2$ must originate from previous
sums and differences, and therefore involve only sine or cosine terms
that may have been transformed by reflection, power, product, and inverse
operators, 
all  of which leave the basic form ${\bf A} {\bf w}(f,g )$ invariant.
 \qed
\end{proof}
\end{theorem}
\begin{corollary}
\label{nice}
An invertible wave number ${\bm \omega}$ may be represented
in terms of multiplicative wave numbers.
\begin{proof}
By the previous theorem, the general form 
of a multiplicative wave number
is the product ${\bf A} {\bm w}(f,g)$ of an amplitude function and a multiplicative wave number.
It follows from Theorem \ref{thm:wnform11} that the general form of the amplitude function
${\bf A}$
involves only functions of expressions of the form
${\bf F}=2 \pi (f {\bm \xi} + g {\bf 1})$. 
But ${\bf F}=-i ln({\bf w}(f,g)$ for some multiplicative wave number
{\bf w}(f,g).\qed
\end{proof}
\end{corollary}
\subsection{Some Elementary Properties of the Invertible Wave Numbers\\}
The elementary properties of invertible wave numbers considered
here include
their field and vector space structure,
their multiplicative norm,
and their transformations.
\begin{theorem}
The set of invertible wave numbers forms 
a field. 
\begin{proof}
It follows from the definition of invertible wave numbers
that the set ${\bf W}_I$ contains 
the set of additive wave numbers. Hence
the identity element ${\bm 0}$ is an invertible wave number.
By the definition of the invertible wave numbers,
the additive inverse $-{\bf A}  {\bf w}(f,g)$ of any invertible
wave number
is in the set ${\bf W}_I$.
Finally it follows from the 
commutativity and associativity of $\oplus$ that 
 ${\bf W}_I$ forms an Abelian additive group under $\oplus$.
It follows from the definition of an invertible wave number
that ${\bf W}_I$ is closed under $\otimes$, which is commutative and associative.
Since the identity element
${\bf w}(m,0)$ is a multiplicative wave number it is also an invertible
wave number and an element of ${\bf W}_I$.
It follows from the definition of the invertible wave numbers
that for any element ${{\bf A}{\bf w}(f,g)}$ 
the multiplicative inverse element
${\bf A}^{-1} \bar{\bf w}(f,g)$ exists.
Hence ${\bf W}_I(\otimes)$ is an Abelian multiplicative group.
Hence the set of invertible wave numbers forms a field ${\mathbb W}_I$.
\end{proof}
\end{theorem}
In addition to its field structure,
the set ${\bf W}_I$ has the structure of a vector space
in which the set of vectors is the set of 
multiplicative wave numbers ${\bf W}_M$ 
together with the element ${\bf 0}$
and 
in which the usual field of scalars is replaced by the
multiplicative group ${\bf F}_I$
of amplitudes specified by
\begin{definition}
\label{def:F}
$
{\bf F}_I\coloneq\big\{ {\bf A}  \mid \forall {\bf w}\ \epsilon\ {\bf W}_M, 
\ {\bf A}  {\bf w} \ \epsilon\ {\bf W}_I.
 \big\}.
$
\end{definition}
\begin{theorem}
\label{th}
${\bf F}_{I}$ 
is a proper multiplicative subgroup of 
${\bf W}_I$.
\begin{proof}
By \ Proposition\ \ref{prop11}, 
${\bf A}\ \epsilon \ {\bf F}_I \Rightarrow {\bf A}\ \epsilon \ {\bf W}_I.$ 
Since the amplitude of ${\bf A}{\bf w}$ is assumed to be in w-free form,
definition \ref{def:F}, prohibits any member of ${\bf W}_M$,
with exception of ${\bf 1}$, from being a factor of any member of ${\bf F}_I$.
Hence ${\bf F}_I \subsetneq {\bf W}_I$ 
and 
${\bf F}_I \cap {\bf W}_M = {\bf 1}$.
If 
${\bf A}, {\bf A}' \ \epsilon \ {\bf F}_I \ then\ \forall\ 
{\bf w},{\bf w}' \ \epsilon \ {\bf W}_M\ one \ has
\\
\ A{\bf w} A'{\bf w}' = \ (A A'){\bf w}{\bf w}'  
\epsilon \ {\bf W}_I\ and\ A  A' \ \epsilon \ {\bf F}_I.
\ By \ Definition \ \ref{def:F}\
A  \ \epsilon \ {\bf F}_I \Rightarrow A^{-1}  \ \epsilon \ {\bf F}_I.
\ Hence \ {\bf F}_I
\ is\ 
a\ multiplicative\  subgroup\ of \  {\bf W}_I.$
\qed
\end{proof}
\end{theorem}
\begin{remark}
${\bf F}_I $ is not an additive subgroup of ${\bf W}_I$.
Note that
if $A, A' \
\epsilon \ {\bf F}_I \ $then $\ \forall \ {\bf w}\ \epsilon \ {\bf W}_I \ $
one\ has $A{\bf w}\oplus\ A'{\bf w} =(A \oplus A'){\bf w} \ \epsilon \ 
{\bf W}_I.$  Setting ${\bf A} = cos({\bf a})$ and 
${\bf A}^{'} ={\bf i} sin({\bf a})$ implies that,
if ${\bf F}_I$ is an additive group, that $e^{i{\bf a}}\ \epsilon\ {\bf F}_I$
which is not permitted by definition. 
\end{remark}
\begin{corollary}
\label{th1}
${\bf W}_I ={\bf F}_I \hat{\otimes}  {\bf W}_M$ 
and ${\bf F}_I \cap {\bf W}_M = {\bf 1}$,
in which $ \hat{\otimes}$ is the set product operator.
\begin{proof}
By Theorem \ref{thm:genform1}, every element of ${\bf W}_I $ has the form
of the product of an amplitude ${\bf A}$ 
and a multiplicative wave number ${\bf w}$.
By definition ${\bf F}_I$ contains every amplitude
${\bf A}$ of ${\bf W}_I$ that does not contain
a member of ${\bf W}_M$ as a factor, apart from ${\bf 1}$, 
and ${\bf W}_I $ contains every 
multiplicative wave number ${\bf w}\ \epsilon\ {\bf W}_M$.  
By the multiplicative closure of ${\bf W}_I $, it must
contain the product of every element
of ${\bf F}_I$ and every element of ${\bf W}_M$ must be an element of
${\bf W}_I $.
By the proof of Theorem \ref{th},
${\bf F}_I \cap {\bf W}_M = {\bf 1}$.
\end{proof}
\end{corollary}
\begin{theorem}
\label{th1}
The multiplicative group ${\bf W}_{M0} ={\bf W}_{M} \cup {\bf 0} $ is a vector space over the  multiplicative group ${\bf F}_I$.
\begin{proof}
The vector ${\bf 0}\ \epsilon\ {\bf W}_{M0}$ by definition.
For all vectors ${\bf w}\ \epsilon {\bf W}_{M0} $ and scalars
${\bf A}\ \epsilon {\bf F}_{A} $, it follows from the properties
of ${\bf W}_I$ that vector addition and
scalar multiplication of vectors are well-defined and satisfy the commutative,
associative, and distributive axioms of a vector space.
\end{proof}
\end{theorem}
The norm 
for invertible wave numbers takes the same form
as the norm for additive wave numbers
$\|A{\bf w}\|=\sqrt{\frac{ \underset {\lambda}{\Sigma}		
a^2({ \xi})}{\lambda}	},	
$ in which $a({ \xi})$ is a component of the amplitude function 
${\bf A({\bm \xi})}$	
and $\lambda$ is the period of the wave number.	
It possesses the multiplicative property.
Employing this form, one may represent
an invertible wave number as having an amplitude that is
a magnitude.
\begin{theorem}
\label{thm:genform}
An invertible wave number is representable
as a product of the norm $\|{\bf A}\|$ of an amplitude
and a multiplicative wave number
\begin{eqnarray}
\label{eq:formgen}
{\bf w}=\|{\bf A}\|e^{i {\bm  \theta}}.
\end{eqnarray}
\begin{proof}
\begin{eqnarray}
{\bf w}&=&{\bf A}{\bf w}(f,g)
=({\bf A}_R \oplus i{\bf A}_I)(\cos (f{\bm  \xi}+g{\bm  1}) 
\oplus i \sin(f{\bm  \xi}+g{\bm  1})\\
&=&({\bm \alpha}  \oplus i{\bm \beta})
\end{eqnarray}
in which ${\bf A}_R$ and ${\bf A}_I$ are the real and imaginary 
components of ${\bf A}$; 
${\bm \alpha}=
{\bf A}_R \cos (f{\bm  \xi}+g{\bm  1})\ominus {\bf A}_I \sin (f{\bm  \xi}+g{\bm  1})$;
and ${\bm \beta}=
{\bf A}_R\sin (f{\bm  \xi}+g{\bm  1}) \oplus{\bf A}_I \cos (f{\bm  \xi}+g{\bm  1})$. 
On noting that
\begin{eqnarray}
\|{\bf w}\|=\sqrt{{\bf A}_R^2 \oplus{\bf A}_I^2}
\end{eqnarray}
one may write
\begin{eqnarray}
{\bf w}=\|{\bf w}\|(\cos ({\bm  \theta}) \oplus i \sin ({\bm  \theta}))
=\|{\bf w}\|e^{i {\bm  \theta}}
\end{eqnarray}
in which ${\bm  \theta}=\cos^{-1} \big(\frac{{\bf A}_R}{\|{\bf w}\|}\big)$.\qed
\end{proof}
\end{theorem}
One may apply the norm in showing how the parameters of
the invertible wave numbers may be extended from rational to real values.
\begin{theorem}
If there is a Cauchy sequence $\rho_n$ of rational parameters  
that converges to
the real parameters $\rho_0$, then 
the sequence
of invertible wave numbers
${\bm \omega}(\rho_n)$ converges to the wave number with real parameters
${\bm \omega}(\rho_0)$.
\begin{proof}
By Theorem \ref{SandD}
\begin{eqnarray}
\|{\bm \omega}(\rho_m)\ominus{\bm \omega}(\rho_n)\| 
&=& 
\| {\bf {2 i}} sin({\bf i} ln\Big(\frac{{\bm \omega}(\rho_n)}
{{\bm \omega}(\rho_m)}^{1/2}\Big)\| \|
\big( {\bm \omega} (\rho_m) \ {\bm \omega}(\rho_n)\big)^{1/2}\|
\nonumber\\
&\rightarrow&    {\bf 0}\ \ as \ n,m\ \rightarrow\ \infty
\end{eqnarray}
\ \ \ \ \ \ \ \ since 
$ {\underset {n,m \Rightarrow \infty}
\lim} ln\Big(\frac{{\bm \omega}(\rho_n)}{{\bm \omega}(\rho_m)}^{1/2}\Big)
=ln({\bf 1})={\bf 0}$ \qed
\end{proof}
\end{theorem}
There are a variety of transformations that may be applied 
to invertible wave numbers, including dilations, rotations, and translations.
\begin{definition}
A dilation ${ D}_\rho$ of an invertible wave  number ${\bm A}{\bf w}(f,g)$ is
\begin{eqnarray}
{ D}_\rho ({\bm A}{\bf w}(f,g))
= {\bm \rho}{\bm A} {\bf w}(f,g),\ \ \ \  {\bm \rho}\ \ \epsilon \ \mathbb{W}_Q 
\end{eqnarray}
\end{definition}
\begin{proposition}
The set of invertible wave numbers is closed under dilation
and ${\bf W}_I=\widehat{D}_{\rho}({\bf W}_I)$
\end{proposition}
\begin{definition}
\label{prop:rot}
A rotation ${ S}_{\xi_0}$ of an invertible wave number ${\bf A}{\bf w}(f,g)$ 
is
\begin{eqnarray}
{ S}_{\xi_0}\big({{ A}(\bm \xi)}{ w}(f,g)(\bm \xi)\big) &=&
{ A}({\bm \xi+\xi_0)}{\bf w}(f,g) (\bm \xi+\xi_0) 
\end{eqnarray} 
when ${\bf A}{\bf w}(f,g)$ 
is represented as a function of ${\bm \xi}$.
\end{definition}
\begin{proposition}
Rotations of wave numbers distribute over 
$\oplus$ and $\otimes$.
\end{proposition}
\begin{proposition}
Rotations of wave numbers form a multiplicative group.
\end{proposition}
%
%
\begin{definition}
A translation ${ T}_{\bf c}$ 
of an invertible wave number $A{\bf w}(f,g)$ is its sum
with a constant wave number ${\bf c}={\bf A}_0{\bf w}(0,g^{\prime}) $
\begin{eqnarray}
{ T}_{\bf c}(A{\bf w}(f,g)) =A {\bf w}(f,g)
 \oplus {\bf A}_0{\bf w}(0,g^{\prime}) 
\end{eqnarray}  
\end{definition}
\begin{proposition}
The translations of invertible wave numbers form a multiplicative group.
\end{proposition}
\subsection{Equations in Invertible Wave Numbers\\}
Since the invertible wave numbers
form a field, the four arithmetic operators may be applied in finding
wave numbers that satisfy equations, as in the case of the simple equation
\begin{example}
The fixed point invertible wave numbers associated with a Mobius
transformation
$M ({\bm \omega}) =( {\bm A} {\bm \omega} \oplus {\bm B})/
({\bm C} {\bm \omega} \oplus{\bm D})$
of an invertible wave number $\bm \omega$,
in which $ {\bm {A,B,C,D}}$ are amplitudes,
are defined by the conditions
\begin{eqnarray}
{\bm \omega}  
&=& {\frac     { {\bm A} {\bm \omega} \oplus{\bm B}}
                 {{\bm C} {\bm \omega} \oplus{\bm D}}},\ \ \ \ \ \ \ 
{{\bm C} {\bm \omega} \oplus{\bm D}} \ne {\bf 0}
\end{eqnarray} 
and hence by the quadratic equation
\begin{eqnarray}
{\bm C} {\bm \omega}^2 \oplus ({\bm D}\ominus{\bm A} ) {\bm \omega}\ominus {\bm B} 
&=& {\bm 0}
\end{eqnarray} 
which may be divided by ${\bf C}$ and solved
for the invertible wave numbers
\begin{eqnarray}
{\bm \omega }
&=& \frac{  -  ({\bm D}\ominus{\bm A} )   \tiny{\begin{array}{c}\oplus\\ \ominus\end{array} }
\sqrt{({\bm D}\ominus{\bm A} ) ^2 
\ominus{\bm 4}  {\bm {BC}} }}
{{\bm 2} {\bf C}}.
\end{eqnarray} 
\end{example}
More generally one may solve equations of the form
$\overset {N} {\underset {j=0} {\oplus} }
{\bf C}_j  {\bf w }_j={\bf 0},\ {\bf C}_j \ne {\bf 0}$ by writing
${\bf C}_j {\bf w}_j = e^{i {\bf F}_j}$,
in which ${\bf F}_j=f_j {\bm \xi} + g_j -i ln {\bf C}_j$,
to obtain  the equation
\begin{eqnarray}
\label{eq:poly}
\overset {N} {\underset {n=0} {\oplus} }
e^{i {\bf F}_j}={\bf 0}
\end{eqnarray}
to which one may immediately apply 
Theorem \ref{thm:wnform11}
to obtain
\begin{eqnarray}
\label{eq:zero}
\overset {N} {\underset {n=0} {\oplus} }
e^{i {\bf F}_j}=
{A}_N({\bf F}_1,...,{\bf F}_N) e^{i \big(\overset {N} {\underset {j=1} \Sigma }\frac{ {\bf F}_j}{N}\big)} ={\bf 0}
\end{eqnarray}
in which
\begin{eqnarray}
A_{N}({\bf F}_1,...,{\bf F}_{N})
&=&
4
\Bigg(
{
\overset {N}
{\underset {m=1}
 \Pi }}A^{\tfrac{1}{2}}_{(N-1)m}
{{\cos\Big({\tfrac{{\overset {N-1}  {  \underset {j\ne m} {\underset {j=1} 
 \Sigma }}{ {\bf F}_j}}-{(N-1)}{\bf F}_m}{2(N-1)}}}-ilnA^{\tfrac{1}{2}}_{(N-1)m}}\Big)\Bigg)^{\tfrac{1}
{N}}\nonumber.
\end{eqnarray}
It is clear that Equation (\ref{eq:zero}) can take on a zero value
for every value of $\xi$ defining the principal sequence
of ${A}_N({\bf F}_1,...,{\bf F}_N) $
iff at least one
of its $N$ factors is zero.
Hence one has
\begin{theorem}
The equation in invertible wave numbers
\begin{eqnarray}
\overset {N} {\underset {n=1} {\oplus} }
{\bf C}_n{\bf w}(f_n,g_n)={\bf 0},\ {\bf C}_j \ne {\bf 0}\ \ 
\end{eqnarray}
has solutions corresponding to the zeros of each of the $N$ conditions
\begin{eqnarray}
{\bf A}^{\tfrac{1}{2}}_{(N-1)m}{{\cos\Big({\tfrac{{\overset {N-1}  {  \underset {j\ne m} {\underset {j=1} 
 \Sigma }}{{\bf F}_j}}-{(N-1)}{\bf F}_m}{2(N-1)}}}-iln{\bf A}^{\tfrac{1}{2}}_{(N-1)m}}\Big)={\bf 0},\ \ \ m=1,N
\end{eqnarray}
\end{theorem}
\begin{example}
It is instructive to consider the theorem for the 
simplest cases of the special case of  $\overset {N} {\underset {j=0} {\oplus} }
{\bf w }_j={\bf 0}$.
For the case $N=2$, one has
\begin{eqnarray}
\label{eq:poly1}
&\ &\overset {2} {\underset {j=1} {\oplus} }
{\bf w }(f_j,g_j)
=
\overset {2} {\underset {j=1} {\oplus} }
e^{  2\pi    i (f_j { \bm \xi} + g_j )}
=
\overset {2} {\underset {j=1} {\oplus} }
e^{i {\bf F}_j} ={\bf 0}\nonumber
\Rightarrow 2cos(\tfrac{  {\bf F}_1-{\bf F}_2 }{2})={\bf 0}
 \nonumber
\end{eqnarray}
and hence
\begin{eqnarray}
2\pi( (f_1-f_2) {\xi} + (g_1 - g_2)) ={(2 \xi -1)}\pi
\ \ \ for \ \xi\ \in\ principal\  sequence
\end{eqnarray}
with solution $f_1=f_2+2$ and $g_1=g_2-1$.
For the case of $N=3$, one has
\begin{eqnarray}
\label{eq:poly1}
&\ &\overset {3} {\underset {j=1} {\oplus} }
 {\bf w }(f_j,g_j)
=
\overset {3} {\underset {j=1} {\oplus} }
e^{      i (f_j { \bm \xi} + g_j )}
=
\overset {3} {\underset {j=1} {\oplus} }
e^{i {\bf F}_j} \nonumber \\
&=&4\big({\bf A}_{21}^{1/2}cos(\tfrac{  {\bf F}_2+{\bf F}_3-2{\bf F}_1  }{2}
-iln{\bf A}_{21}^{1/2})
({\bf A}_{22}^{1/2}cos(\tfrac{  {\bf F}_3+{\bf F}_1-2{\bf F}_2  }{2}
-iln{\bf A}_{22}^{1/2}) \nonumber \\
&\ &\ \ ({\bf A}_{23}^{1/2}cos(\tfrac{  {\bf F}_1+{\bf F}_2-2{\bf F}_3  }{2}
-iln{\bf A}_{23}^{1/2}) 
\end{eqnarray}
in which 
${\bf A}_{21}=
2cos(\tfrac{f_2-f_3+g_2-g_3}{2}-iln (\tfrac{{\bf A}_2} {{\bf A}_3})^{1/2}),
{\bf A}_{22}=2cos(\tfrac{f_1-f_3+g_1-g_3}{2}-iln (\tfrac{{\bf A}_1} {{\bf A}_3})^{1/2})$, and
${\bf A}_{23}=2cos(\tfrac{f_1-f_2+g_1-g_2}{2}-iln (\tfrac{{\bf A}_1} {{\bf A}_2})^{1/2})$.
One first notes that each of the three products in Equation \ref{eq:poly1} 
is zero if
${\bf A}_{2j}=0$ for the corresponding value of $j$.
In terms of the first product ${\bf A}_{21}$, for example, this implies
that $ln {\bf A}_2/{\bf A}_3=0$, and hence that ${\bf A}_2={\bf A}_3$,
and that $cos((f_2-f_3)\xi +(g_2-g_3))=0$, whence $f_3=f_2-2 \pi$ and 
$g_3=g_2-\pi$. There is another solution corresponding to
${\bf A}_{2j}=1$, which implies that ${\bf A}_2={\bf A}_3$ and that 
$\tfrac{f_2+f_3-2f_1)\xi+(g_2+g_3-2g_1)}{4}-i ln {\bf A}_2{\bf A}_3/{\bf A}_1^2=0$, which in turn implies that ${\bf A}_1={\bf A}_2={\bf A}_3$
and that $f_2+f_3-2f_1=4\pi$ and that $g_2+g_3-2g_1=2\pi$.
\end{example}
\subsection{Orthonormal Bases and
Particulate, Integral, and Prime Wave Numbers\\}
One may extend the $n$ orthogonal bases ${\bm u}(\tfrac{1}{n},j),j=1,n$ 
of additive wave numbers
defined in Section \ref{sec:bases} to bases of invertible wave numbers
that are analogous 
to the standard bases of Euclidean geometry.
\begin{theorem}
\label{cor:basis}
For every $n\ \epsilon\ 
\mathbb{N}^+$,
there exist $n$ invertible wave numbers 
\begin{eqnarray}
\Big\{
{\bm e}(\tfrac{1}{n},\xi)=
(0,...,0,1, 0,...,0 )\mid \xi=1,n \Big\},
\end{eqnarray}
in which the 1 appears in the $\xi$th place,
that form an orthonormal  basis for representing
an invertible wave number ${\bf A} {\bf w}(f,g)$.
\begin{proof}
Using Equation \ref{eq:basis} define \begin{eqnarray}
{\bm e}(\tfrac{1}{n},j)=\tfrac{1}{\bf n} {\bm u}(\tfrac{1}{n},j),j=1,n
\end{eqnarray}
from which it follows that the jth term of the principal sequence
is
\begin{eqnarray}_{1}{\bm e}(\tfrac{1}{n},j)=(0,...,1,...,0)\end{eqnarray}
with the $1$ occurring in the jth position.
Clearly the $n$ numbers form an orthonormal set 
of invertible wave number under $\otimes$
and it follows from Definition \ref{def:invertible} for the invertible wave numbers
that the expression
${\underset {j \epsilon S}
\oplus}{\bm e}(\tfrac{1}{n},j)({\bf A}{\bf w}(f,g))$, in which 
${\bf A}{\bf w}(f,g)$ an invertible wave number
of period $n$ and $S$ is a subset of the first $n$ natural numbers, 
is an invertible wave number of period $n$ with
$
\overset {n} {\underset {j=1}
\oplus}
{\bf A}{\bf w}(f,g)
{\bf e}(\tfrac{1}{n},\xi)
={\bf A} {\bf w}(f,g).
$
\qed
\end{proof}
\end{theorem}
The wave numbers of the basis $\big\{{\bm e}(\tfrac{1}{n},j), j=1,n\big\}$ 
may be viewed as operators under $\otimes$ 
that select
the phases that comprise a wave number.
One may therefore employ these operators
in constructing wave numbers
that contain arbitrary subsets of the phases of an invertible wave number
${\bf A}{\bf w}(f,g)$
and zeros in place of the non-included phases.
Three interesting examples of such classes
are the integral, particulate,
and prime wave numbers.
\subsubsection{The Integral Wave Numbers}
Integral wave numbers are sequences of numbers whose phases
are the cumulative sums of the phases of a wave number.
The integral wave numbers may be constructed by applying
the orthonormal basis to previously-defined wave numbers in a manner 
specified in
\begin{definition}
The principal value of the integral wave number
${ \smallint}{\bf A}{\bf w}(f,g)$
of the invertible wave number 
${\bf A}{\bf w}(f,g)$ of period $n$ is
\begin{eqnarray}
{_1}{ \smallint}{\bf A}{\bf w}(f,g)
&=&
\big\{
\overset {k} {\underset {j=1}\oplus}
{\bf A}{\bf w}(f,g)
{\bf e}(\tfrac{1}{n},j) \mid  k=1,n \big\}
\end{eqnarray}
\end{definition}
from which follows
\begin{theorem}
The integral $\smallint {\bf A}{\bf w}$ of an invertible wave number is an invertible wave number.
\begin{proof}
Writing the $j$th phase of ${_1}{\bf A}{\bf w}$ as ${ a(j)}e^{2 \pi i (\tfrac{m}{n}j+g)}$,
the $k$th phase of $\smallint {\bf A}{\bf w}$ may be written
\begin{eqnarray}
I_k
&=&%
\overset {k} {\underset {j=1}\Sigma}
{ a(j)}e^{2 \pi i (\tfrac{m}{n}j+g)}=
\overset {k} {\underset {j=1}\Sigma}
e^{i(-iln a(j) + 2\pi (\tfrac{m}{n}j+g))}
\end{eqnarray}
to which one may apply Theorem \ref{thm:wnform11} to obtain
\begin{eqnarray}
I_k
&=&
\hat{A}(k)
{e^{i2\pi{{{\overset {k}  { {\underset {j=1} 
 \Sigma }}\big(\tfrac{m}{n}}j+g\big)}}}}
=
\hat{A}(k)
e^{i2\pi\big(
\tfrac{m}{2n}(k+1)+g\big)
}
=
\hat{\hat{A}}(k)
e^{i2\pi\big(
\tfrac{m}{2n}(k)+g\big)
}\\
&=&
\hat{\hat{A}}(k)
{\bf w}(\tfrac{m}{2n},g)
\end{eqnarray}
which is an invertible wave number that has been 
scaled, translated, and rotated. \qed
\end{proof}
\end{theorem}
Integral wave numbers formed from 
the simplest multiplicative wave numbers 
$
{_1}{\smallint}({\bf w}(\tfrac{1}{n},0))
=
\big\{
\overset {k} {\underset {j=1}\oplus}
{\bf w}(\tfrac{1}{n},0){\bf e}(\tfrac{1}{n},j\big), k=1,n \big\}
$
take the form noted in
\begin{theorem}
\label{recursngon}
The integral wave number $
{\smallint}({\bf w}(\tfrac{1}{n},0))$
has a geometric representation as a regular n-gon with a vertex at the origin.
\begin{proof}
Interpreting the wave number ${\bf w}(\frac{1}{n},0)$ as a regular n-gon,
one may view the construction of the integral wave number
as the successive addition of unit vectors 
$\overset {k} {\underset{ j=1}\Sigma} e^{i2\pi \tfrac{j}{n}}, k=1,n$ whose arguments
differ by the same angle $2\pi/n$ and form the $n$ exterior angles 
of a geometric form. Since the sum of these angles is $2\pi$
and since the vectors are of the same unit length, they must form
a regular n-gon. Since the first unit vector $e^{i2\pi \tfrac{j}{n}}$
originates at the origin, the last one $e^{i2\pi \tfrac{n}{n}}$ 
must terminate there, and the origin is a vertex of the n-gon. \qed
\end{proof}
\end{theorem}
One may iterate this process
and recursively construct 
countable sequences of multiplicative wave numbers
by choosing the norm $\|{\bf w}(\tfrac{1}{n},0)\|$
of the multiplicative wave number 
at one iteration to serve as the
length of the edge of the multiplicative wave number
in the following iteration, a condition that is automatically
satisfied in the first iteration specified by Theorem \ref{recursngon}.
\begin{theorem}
\label{thm:circular1}
For every multiplicative wave number ${\bf w}( \tfrac{1}{n},0)$, 
$n\ \epsilon\ \mathbb{N}$
there exists a countable, recursive sequence of integral wave numbers
whose elements have the principal sequence
\begin{eqnarray}
{_1}{ I^t}({\bf w}( \tfrac{1}{n},0))
=
\overset {k} {\underset {\xi=1}\oplus}
{_1}{ I^{t-1}}\big(({\bf w}( \tfrac{1}{n},0){\bf e}(\tfrac{1}{n},\xi)\big), \ \ \  t\ \epsilon\ \mathbb N,
\end{eqnarray}
with initial value $I^0({\bf w}(\tfrac{1}{n},0))={\bf w}(\tfrac{1}{n},0)$,
and norm 
$\|{ I^t}({\bf w}( \tfrac{1}{n},0)\|=
\tfrac{\|{ I^{t-1}}({\bf w}( \tfrac{1}{n},0))\|}{2\sin\big(\tfrac{\pi }{n}\big)}$
\begin{proof}
One may view the proof of Theorem \ref{recursngon} as describing the
initial step of the recursive process. In each following iteration an analogous
procedure is applied, but to a regular n-gon whose edges have a length equal 
to the magnitude
of the $n$ defining vectors of the previous n-gon
and defining vectors that have been rotated and translated to a new origin.
This process is described by the first order,
linear recurrence relation
\begin{eqnarray}
r_0=1;\  e_t=r_{t-1}; \ e_t=2 \sin (\tfrac{\pi}{n})r_t
\end{eqnarray}
in which $r_t$ is the norm of the defining vectors 
at iteration $t$ and $e_t$ is the length
of an edge at time $t$.
Hence at any iteration $t$, the process generates 
a regular n-gon whose edges have length $e_t$
and subtend an angle $2\pi/n$ at the center of the n-gon.
The $n$th  element of the sequence
$\underset {j=1} {\overset {n} \oplus }{\bm u}(\tfrac{1}{n}, j)=0$ since
the $n$ $n$th roots of unity sum to zero, 
while the $(n-1)$st element 
$\underset {j=1} {\overset {n-1} \oplus }{\bm u}(\tfrac{1}{n}, j)=-1$ 
since ${\bm w}(\tfrac{1}{n},n)=1$.
\end{proof}
\end{theorem}
\begin{corollary}
For $n=6$, the norm of the wave numbers at each iteration is invariant;
for $n<6$ the norm decreases monotonically at each iteration;
and increases monotonically for $n>6$.
\end{corollary}
\begin{remark}
The $n$ translates of a wave number at any iteration of Theorem \ref{thm:circular1}
may be formed into another basis for the invertible wave numbers
by forming the $n$ products of $n-1$ of the translates,
as in Theorem \ref{thm:basis}
\end{remark}
The property of integral wave numbers that the last component
$\underset {j=1} {\overset {n} \oplus }{\bm u}(\tfrac{1}{n}, j)=0$
is the equivalent for wave numbers
of the fact that the $n$ $n$th roots of unity sum to zero
and
leads to the
\begin{definition}
An additive wave number is zero-sum iff the sum 
of its phases is zero.
\end{definition}
\begin{proposition}
\label{theor:wsum0}
Multiplicative wave numbers
of $period\geq 2$
are zero-sum.
\begin{proof}
Let $c(k)=e^{2 \pi i\big( \tfrac{k}{n} \big)}$ for $k=1,n$ then
\begin{eqnarray}
\overset {n} {\underset {k=1} {\Sigma} }c(k)=
\overset {n} {\underset {k=1} {\Sigma} }c(1)^k=
\frac{c(1)-c(1)^{k+1}}{1-c(1)}=\frac{0}{1-c(1)}=0
\end{eqnarray}
since $c(1)^{k+1}=c(1)$.
Hence for the multiplicative wave number ${\bf w}(\tfrac{m}{n},\tfrac{p}{q})$,
the sum of the phases is
\begin{eqnarray}
\underset {k=1} {\overset{nq} {   \oplus} }
{\bf w}(\tfrac{m}{n},\tfrac{p}{q})
{\bf w}(0,\tfrac{k}{nq})=
{\bf w}(\tfrac{m}{n},\tfrac{p}{q})
\underset {k=1} {\overset{nq} {   \oplus} }
{\bf w}(0,\tfrac{k}{nq})
=
{\bf w}(\tfrac{m}{n},\tfrac{p}{q})\underset {k=1} {\overset{nq} {   \Sigma} }
e^{2 \pi i\big( \tfrac{k}{nq} \big)}
=0
\end{eqnarray}
\qed
\end{proof}
\end{proposition}
\begin{proposition}
\label{theor:wsum2}
The sum and product of zero-sum wave numbers
are zero-sum wave numbers.
\begin{proof}
For zero-sum sequences ${\bf A_1}{\bf w}(f_1,g_1)$ and 
${\bf A_2}{\bf w}(f_2,g_1)$, one has
\begin{eqnarray}
S_1=\overset {\lambda_1} {\underset {\xi=1} {\oplus} }
{\bf A_1}{\bf w}(f_1,g_1)=0;\ \ 
S_2=\overset {\lambda_2} {\underset {\xi=1} {\oplus} }
{\bf A_2}{\bf w}(f_2,g_2)
=0
\end{eqnarray}
in which $\lambda_k$ is the period of ${\bf A_k}{\bf w}(f_k,g_k)$.
By the definition of $\oplus$ and $\otimes$ for two sequences of period 
$\lambda_1$ and $\lambda_2$, it follows that 
\begin{eqnarray}
\big({\bf A_1}{\bf w}(f_1,g_1)\oplus{\bf A_2}{\bf w}(f_2,g_2)\big)&=&\lambda_2 S_1 + \lambda_1 S_2=0\\
\big({\bf A_1}{\bf w}(f_1,g_1)\otimes{\bf A_2}{\bf w}(f_2,g_2)\big)&=&S_1 S_2=0
\end{eqnarray}
\qed
\end{proof}
\end{proposition}
\subsubsection{The Particulate Wave Numbers\\}
One may construct a basis for the wave numbers
from the $n$ bases for the $n$th roots of unity
that allows one to construct a wave number that takes the form
of a single, isolated value,
and is therefore not periodic.
In proving this result, it is of value to introduce
the re-and co-numbers of a wave number, which represent
two special cases of choosing which phases 
to select in defining new wave numbers from existing wave numbers.

In terms of basis wave
numbers ${\bf e}(\tfrac{1}{n},j)$, one has
\begin{definition}
\label{co-number}
The principal sequences of the co- and re-numbers of a basis wave number
${\bf e}(\tfrac{1}{n},j)$ are, respectively
\begin{eqnarray}
_1\overset {\bm \ast}{\bf e}(\tfrac{1}{n})
&=&_1\Big(\overset {n-1} {\underset {\xi=1}\oplus}
{\bm e}(\tfrac{1}{n},\xi) \Big)\sim \big( 1,...1,0\big)\\
_1\overset {\bm \circ}{\bf e}(\tfrac{1}{n})
&=&
_1{\bm e}(\tfrac{1}{n},n)\sim \big( 0,...,0,1\big)\
\end{eqnarray}
\end{definition}
\begin{proposition}
\label{prop:binary0}
$\ \ \ \overset {\bm \ast}{\bf e}(\tfrac{1}{n})
\oplus \overset {\bm \circ}{\bf e}(\tfrac{1}{n})
={\bf 1};\ \ \ \ \ \ \overset {\bm \ast}{\bf e}(\tfrac{1}{n})
\otimes \overset {\bm \circ}{\bf e}(\tfrac{1}{n})
={\bf 0}.$
\end{proposition}
\begin{definition}
The symmetric particulate wave numbers are
\begin{eqnarray}
{\bf P}_n(1)
&=&\overset {\bm \circ}{\bf e}(\tfrac{1}{n})
\big(\overset {}{\underset {j>n}\otimes}
\overset {\bm \ast}{\bf e}(\tfrac{1}{j})\big)
 \  \ \ \   -\infty < n < \infty,\ n\ne 0.\\
{\bf P}_n(m)
&=&\overset {m} {\underset {k=1}\oplus}
{\bf P}_n(1)
 \  \ \ \  \ \ \ \ \ \ \ \ \    0 < m < \infty.
\end{eqnarray}
\end{definition}
\begin{proposition}
The symmetric particulate wave number ${\bf P}_n(1)$ is a sequence over ${\bm \xi}$ that
takes the value $1$ at $ \xi=\pm n$ and zero elsewhere.
\begin{proof}
The wave number $\overset {\bm \circ}{\bf e}(\tfrac{1}{n})$ has 
the value $1$ at $\xi=n$ and at all multiples of $n$
and zeros elsewhere, while the wave number 
$\overset {}{\underset {j>n}\otimes}
\overset {\bm \ast}{\bf e}(\tfrac{1}{j})$ has the value $1$ for all indices less
than or equal to $n$, and zero for all indices greater than $n$.
Hence their product satisfies the statement of the proposition.
\end{proof}
\end{proposition}
\begin{proposition}
There exist particulate integer numbers ${\bf P}^{+/-}_{n}(m)$ 
that take 
the value $m\ \epsilon\ \mathbb{Z}$ at 
$\xi=+/- n$ respectively and zero elsewhere and
particulate rational numbers ${\bf P}_{n}^{+/-}(q)$
that take the value $q\ \epsilon\ \mathbb{Q}$ at $\xi=\pm n$
and zero elsewhere.
\begin{proof}
\begin{eqnarray}
{\bf P}^+_{n}(m)&=&\tfrac{-e^{i 2 \pi 2 n }{\bf P}_{n}(m)\otimes{\bf P}_{n}(m)}{\bf m};
\ \ 
{\bf P}^-_{n}(m)=\tfrac{e^{-i 2 \pi 2 n }{\bf P}^-_{n}(m)
\otimes{\bf P}^+_{n}(m)
}{\bf m}, \  m\ \epsilon\ \mathbb{N} 
\\
{\bf P}^+_{n}(q)
&=&\tfrac{-e^{i 2 \pi 2 n }{\bf P}^-_{n}(m)\otimes{\bf P}^+_{n}(m)}{{\bf m}{\bf n}};
\ \ 
{\bf P}^-_{n}(q)=\tfrac{e^{-i 2 \pi 2 n }{\bf P}^+_{n}(m)
\otimes{\bf P}^-_{n}(m)}{{\bf m}{\bf n}},
 \ \ q=\tfrac{m}{n}\ \epsilon\ \mathbb{Q} \nonumber
\end{eqnarray}
\end{proof}
\end{proposition}
\begin{remark}
Since one may interpret the wave numbers as numbers
in which all phases hold simultaneously, one may interpret particulate integer
and rational wave numbers
 as wave numbers in which all ambiguity about phase value is resolved
and, as a consequence, have no finite periodicity.
\end{remark}
\subsubsection{The Re-and Co-Numbers and the Period-prime Wave Numbers\\}
In order to specify the period-prime wave
numbers, it is of value to define the re-and co-numbers of 
a multiplicative wave number ${\bf w}(f,g)$:
\begin{definition}
\label{co-number}
The co-number of a simple multiplicative wave number
${\bf w}(\tfrac{m}{n}, \tfrac{p}{n})$ of period $n$ is
\begin{eqnarray}
\overset {\bm \ast}{\bf w}(\tfrac{m}{n}, \tfrac{p}{n}))
&=&\overset {n-1} {\underset {\xi=1}\oplus}
{\bm e}(\tfrac{1}{n},\xi){\bf w}(\tfrac{m}{n}, \tfrac{p}{n}))
\end{eqnarray}
while its re-number 
is
\begin{eqnarray}
\overset {\bm \circ}{\bf w}(\tfrac{m}{n}, \tfrac{p}{n}))
&=&
{\bm e}(\tfrac{1}{n},n){\bf w}(\tfrac{m}{n}, \tfrac{p}{n}))
\end{eqnarray}
\end{definition}
\begin{proposition}
\label{prop:binary}
The co- and re-numbers $
\overset { \ast}{\bf w}(f,g))
$
and $
\overset { \circ}{\bf w}(f,g))
$
partition a multiplicative wave number into
a wave number comprising
its complex and negative terms and a wave number
comprising its positive real term.
\begin{proof} From Definition \ref{co-number}, it is clear that
\begin{eqnarray}
\overset {\bm \ast}{\bf w}(f,g)\oplus  
\overset {\bm \circ }{\bf w}(f,g)
&=&{\bf w}(f,g) \nonumber \\
\overset {\bm \ast}{\bf w}(f,g)\otimes  
\overset {\bm \circ }{\bf w}(f,g)&=&{\bm 0}. 
\end{eqnarray} 
\end{proof}
\end{proposition}
The co-numbers of a 
wave number may be interpreted as specifying values
that the wave number does not take.
This may be exemplified in the problem of identifying arbitrarily large
sets of prime numbers, whose usual
negative definition
(see, for example, Friedlander and Iwaniec (2010)\nocite{FI10})
views them
as natural numbers that are not divisible by preceding prime numbers.
This leads to the problem of representing countable sequences
of numbers 
that are defined in terms of values that they do not possess,
which suggests an application of the co-numbers.

It is convenient to represent the connection between 
the co-numbers and the prime numbers
in terms of 
a circular product operator $\odot$
\begin{definition}
\label{def:circ_prod}
The circular product $\odot$ of multiplicative  wave numbers is
\begin{eqnarray}
{\bf w}  ( \tfrac{  m_1 }{n_1},g_1) \odot {\bf w}
 ( \tfrac{  m_2 }{n_2},g_2)
&=&{\it R}_A
\left({\bf w}  ( \tfrac{  m_1 }{n_1},g_1) \otimes{\bf w}
 ( \tfrac{  m_2 }{n_2},g_2) \right)\nonumber \\
&=&{\bf w} 
\big(\tfrac{  { 1}    }
                   {   n_1n_2    },\tfrac{g_1+g_2}{A}\big)
\end{eqnarray}
in which $A={ m_1n_2+m_2n_1}$ and $R_A$ is the root operator.
\end{definition}
and a cumulative circular product of wave co-numbers
\begin{definition}
The cumulative circular product of wave co-numbers is
\begin{eqnarray}
\label{eq:cumseq}
{\overset \ast{\bm P}}(N)
=
{{\overset {N}  { 
{\underset {n=2} 
 \odot}}}} 
\overset {\bm \ast}{\bf w}(\tfrac{1}{n},0), 
\hspace{.2in} N\ \epsilon\ {\mathbb N}, \ \ 
\end{eqnarray}
\end{definition}
One notes that the zeros of the wave number ${\overset \ast{\bm P}}(N)$ 
occur at, and only at, phases $\xi$ that are
multiples of each of the values $n=2,N$.

Employing these definitions, together with the 
notation $\pi(N)$ for the number of primes occurring in the interval $2,...,N$
one notes the useful property
\begin{lemma}
\label{lem:cumprod}
\begin{eqnarray}
{\overset \ast{\bm P}}(N)
=
{\overset {N}  {  {\underset {n=2}  \odot}}} 
\overset {\bm \ast}{\bf w}(\tfrac{1}{n},0)
=
{\overset {\pi(N)}  { 
{\underset {k=1} 
 \odot}}} 
\overset {\bm \ast}{\bf w}(\tfrac{1}{p_k},0)
\end{eqnarray}
in which $p_k$ denotes the kth prime number.
\begin{proof}
If a co-number for a non-prime integer occurs in the product
$\overset {N}  {\underset {n=1}   \odot}  
\overset {\bm \ast}{\bf w}(\tfrac{1}{n},0)$
it must be preceded in the product by co-numbers for
all of its
prime factors.
The zeros of the sequence of the non-prime co-number
relative to any of its factors
must be a proper subset of the zeros of that factor
by Definition \ref{co-number} of the wave co-numbers.
Hence the product sequence is unaffected by the sequence of any non-prime 
co-number, which may therefore be omitted from the sequence,
leaving a product of only prime co-numbers. \qed
\end{proof}
\end{lemma}
\begin{corollary}
\label{cor:cumprod}
The phase numbers associated with the zeros of the sequence
${\overset \ast{\bm P}}(N)$ correspond to phase numbers 
that are multiples of the prime numbers $\big\{2,...,p_{{\pi}(N)}\big\}$.
\end{corollary}
Since ${\overset \ast{\bm P}}(N)$
may be viewed as a product of prime co-numbers, as well as of
a product of the natural co-numbers, and since one is interested 
in which values of its phase numbers are prime,
it is convenient to represent it as
a function of the phase $\xi$ and denote it with
\begin{eqnarray}
\label{eqn:recursion1}
{\overset \ast{\bm U}_{\pi_{(N)}}}(\xi)
=
\overset { \pi_{(N)} }  {  \underset {k=2}  \odot} 
\overset {\bm \ast}{\bf w}(\tfrac{1}{p_k},0)(\xi)
= \overset {\bm \ast } {\bm U}_{\pi(N)-1}(\xi)
\  {   \odot} \ 
\overset {\bm \ast } {\bf {w}}(\tfrac{1}{p_{\pi(N)}},0)(\xi),
\end{eqnarray}
so emphasizing both its recursive structure
and the dependence of its values on the phase number $\xi$.

One may view the process of recursively constructing the cumulative
product as equivalent to sequentially finding the prime values of
the phase numbers and hence the prime wave numbers.
The phases with values of zero correspond to non-prime phases,
while the phases with non-zero values are candidates
to be prime phases.
It follows that one may employ
circular products of prime numbers
in identifying both the prime wave numbers and the prime natural numbers.
These ideas are made precise in 
\begin{theorem}
\label{thm:big}
If $ p_1,  ..., p_{N+1}$ are the first $N+1$ prime phases,
then the phases in the range $p_{N+1} \leq \xi < p^2_{N+1}$
that are associated with the non-zero terms of the circular product
of the first $N$ prime co-numbers
\begin{eqnarray}
 \overset {\bm \ast }{\bm {\bm U}}_{p_N}
= \underset {n=1} {\overset{ {N}}  {   \odot} }
 \overset {\bm \ast } {\bf {w}}(\tfrac{1}{p_n},0)
\end{eqnarray}
are, together with $ p_1,  ..., p_{N}  $, 
all of the prime phases 
less than $p^2_{N+1}$.
\begin{proof}
Let 
$  \overset {\bm \ast }{\bm {\bm U} }_{p_N}(\xi)$ be 
be values defined on the phases of $  \overset {\bm \ast }{\bm {\bm U}}_{p_N}$.
By Corollary \ref{cor:cumprod},
$\overset {\bm \ast }{\bf w} (\tfrac{1}{p_N},0)(\xi)=0$ 
iff the phase $\xi$ is a multiple of any of the prime numbers 
$p_1,...,p_N$.
The phase $\xi=p^2_{N+1}$ is not prime and is associated with
a non-zero value since it does not satisfy
this condition.
Any phase smaller than $\xi=p^2_{N+1}$ may be represented
by the Fundamental Theorem of Arithmetic
\begin{eqnarray} 
\xi= p_1^{a_1}\cdot p_2^{a_2} \cdot\cdot\cdot p_N^{a_N} 
\cdot {p_{N+1}}^{a_{N+1}}\cdot\cdot\cdot( \phi_{p_{N+1}^2})^{a_m}
\end{eqnarray}
for some $m$,
in which $ \phi_{p_{N+1}^2}$
represents the largest prime phase that is less than ${p_{N+1}^2}$.
Since the function values of phases that contain 
any prime factors less than $p_{N+1}$
are zero
and since permissible phases are less than $p^2_{N+1}$, it follows 
that phases associated with non-zero function values
must have representations
\begin{eqnarray} 
 \xi= p^{s_1}_{N+1}\cdot\cdot\cdot 
( \phi_{p_{N+1}^2})^{s_m}
\end{eqnarray}
for some $m$ in which at most one of the $s_i$'s can take on the value one.
It follows that any non-zero phase that is less than $p_{m+1}^2$
must be a prime number in the range 
${p}_{N+1},..., \phi_{p_{N+1}^2}$. 
\qed
\end{proof}
\end{theorem}
Equation \ref{eqn:recursion1} implies a recursive procedure
for determining the primes
which are defined entirely in terms of the initial condition
$ \overset {\bm \ast } {\bf {w}}(\tfrac{1}{p_2},0)$.
Furthermore, this procedure leads to
the identification of every prime by
\begin{corollary}
The recursive procedure that is defined by
applying Theorem \ref{thm:big} sequentially for $N=1,2,3,...$
does not terminate.
\begin{proof}
Theorem \ref{thm:big} implies that in order to identify prime phases
at step $N+1$, it is necessary at step $N$ to  identify $p_{N+1}$
in order to construct $ \overset {\bm \ast } {\bf {w}}(\tfrac{1}{p_{N+1}},0)$ 
and to identify $p_{N+2}$ in order to ensure the primality
of the numbers generated by $ \overset {\bm \ast } {\bf {w}}(\tfrac{1}{p_{N+1}},0)$.
As shown above, this is the case for the first three steps $N\leq3$.
Assume it to be true for the $N$th step and that $p_{N+1}$ 
and $p_{N+2}$ are known.
Theorem \ref{thm:big} also implies that all prime phases $\xi< p_{N+2}^2$
are identified by 
$ \overset {\bm \ast } {\bf {w}}(\tfrac{1}{p_{N+1}},0)$,
and hence $p_{N+2}$ is identified
at the $(N+1)$st iteration.
By Bertrand's Theorem, $p_{N+3} \leq 2 p_{N+2}$,
and since $2 p_{N+2} \leq  p^2_{N+2}$ for $p_{N+2}>2$,
$p_{N+3}$ is determined. \qed
\end{proof}
\end{corollary}
Theorem \ref{thm:big} implies that one may identify
increasingly-large sets of maximum size of prime numbers
in an iterative process by which, at the $N$th step,
one employs the maximum number
of correct primes
determined in the previous steps to compute
the following largest, error-free set of primes.
In particular, one has
\begin{corollary}
The largest prime phase
in the set of prime phases that may be correctly identified at each application of
Theorem \ref{thm:big}  is 
\begin{eqnarray}
\lfloor 7 \rfloor_p, \ \ 
\lfloor 7^2 \rfloor_p,\ \ 
\lfloor \lfloor 7^2 \rfloor_p^2 \rfloor_p ,\ \
\lfloor \lfloor \lfloor7^2 \rfloor^2_p \rfloor^2_p \rfloor_p,\ \  
\lfloor\lfloor\lfloor \lfloor7^2 \rfloor^2_p \rfloor^2_p\rfloor^2_p \rfloor_p ,\ \ 
\lfloor\lfloor \lfloor\lfloor7^2 \rfloor^2_p \rfloor^2_p\rfloor^2_p \rfloor^2_p\ ,...
\end{eqnarray}
\begin{proof}
In the first step at $N=1$, 
${\overset {\bm \ast } {\bm w }}_1
={\overset {\bm \ast }{\bm U }}_2$ and
it follows from equation (21) that the largest prime phase 
found is $7$. From Theorem \ref{thm:big} it follows that
the largest correctly identifiable prime phase
at each iteration for $N\geq 2$ is the largest prime phase
that is less than the square of the next known prime phase.
Hence for $N=2$, the next prime phase is defined to be $7= \lfloor 7 \rfloor_p$
and largest prime phase that is correctly identifiable is
$\lfloor 7^2 \rfloor_p$,
which is the next prime phase at the $N=3$rd
iteration. 
Hence the largest prime phase identifiable
at $N=3$ is $\lfloor \lfloor 7^2 \rfloor_p^2 \rfloor_p $,
which is the next prime phase at the step $N=4$.
The same argument applies at each iteration $N \geq 4$. \qed
\end{proof}
\end{corollary}
The length of the
domain of phases over which all prime phases are correctly
inferred
therefore increases at each iteration as approximately
the square of the length of the previous domain,
and
is representable as a function of $7$
at each iteration. Recalling that $7$ is the largest
of the prime numbers derived at the first iteration
and that the derivation of the prime numbers is an initial
value problem, this reflects the dependence of the solution
on the initial conditions.

When applied to the results of Corollary $5.5$,
the Gauss/Legendre approximation 
$\pi(N)\approx N/ln(N)$
to the prime counting function, together with 
the approximation
$\lfloor N \rfloor_p \approx N$,  implies
\begin{corollary}
The maximum number of prime phases identifiable at each 
iteration of Theorem $ 5.3$ is approximately
${7^{2(N-1)}}/{(2(N-1)ln7)}$ for $N>1$.
\end{corollary}
It is to be noted that one may employ the re- numbers
$\overset{\circ}{\bf w}(\tfrac{1}{p},0)$ instead of the co-numbers
$ {\overset {\bm \ast }{\bf  w}(\tfrac{1}{p},0)}$ in identifying
prime phases and one has, for example,
\begin{corollary}
If $ p_1,  ...,, p_{N+1}$ are the first $N+1$ prime phases
then $p_N < k <p^2_{N+1}$ is
prime iff
$ 
\underset {j=1} {\overset{ N} \sum}
\overset {\bm \circ } 
{\bf {w}}(j) (\xi)=0.
$
\begin{proof}
By Theorem \ref{thm:big}, the definition of $\overset {\bm \circ } 
{\bm {w}}(j)$, and the definition of the arithmetic summation operator. 
\qed
\end{proof}
\end{corollary}


\newpage

\begin{thebibliography}{00}

\bibitem{FI10}
Friedlander, J. and H. Iwaniec (2010).
"Opera de Cribo."  
American Mathematical Society, Colloquium Publications, Volume 57.

\bibitem{G96}
Goldrei, D.C. (1996).
"Classic Set Theory."  
Chapman and Hall.
296 pages

\bibitem{L92}
Lang, Serge (1992).
Algebra, 3rd edition.
Addison Wesley Publishing Co.
906 pages

\bibitem{MV12}
 Hugh L. Montgomery, H.L and R.C. Vaughan (2012)
Multiplicative Number Theory I: Classical Theory 
(Cambridge Studies in Advanced Mathematics, Series Number 97.

\end{thebibliography}
\end{document}